\documentclass[letterpaper]{amsart}

\pdfoutput=1

\newtheorem{theorem}{Theorem}[section]
\newtheorem{lemma}[theorem]{Lemma}
\newtheorem{proposition}[theorem]{Proposition}
\newtheorem{corollary}[theorem]{Corollary}
\usepackage[top=1.0in, left=1.0in, right=1.0in]{geometry}
\usepackage{graphicx}
\usepackage[all,cmtip]{xy}
\usepackage{xypic}

\theoremstyle{definition}
\newtheorem{definition}[theorem]{Definition}
\newtheorem{example}[theorem]{Example}

\newtheorem{remarks}[theorem]{Remarks}

\theoremstyle{remark}
\newtheorem{remark}[theorem]{Remark}
\usepackage{hyperref}
\usepackage{graphicx}
\usepackage{color}
\usepackage{amsmath}
\usepackage{amsfonts}
\usepackage{euler}
\usepackage{amssymb}
\usepackage{epsfig}
\usepackage{rotating}
\usepackage{graphics}
\newcommand{\be}{\begin{equation}}

\newcommand{\ee}{\end{equation}}
\newcommand{\ben}{\begin{equation*}}
\newcommand{\een}{\end{equation*}}

\newcommand{\wh}[1]{\widehat{#1}}

\newcommand{\om}{\omega}

\newcommand{\bal}{\begin{aligned}}
\newcommand{\eal}{\end{aligned}}

\newcommand{\f}{\varphi}

\newcommand{\dz}{\wedge}

\newcommand{\C}{\mathcal{C}}

\newcommand{\ba}{\begin{array}}

\newcommand{\ea}{\end{array}}

\newcommand{\beq}{\begin{eqnarray}}

\newcommand{\eeq}{\end{eqnarray}}

\newtheorem{lm}{lemma}

\newtheorem{thee}{theorem}

\newtheorem{proo}{proposition}

\newtheorem{co}{corollary}

\newtheorem{rem}{remark}

\newtheorem{deff}{definition}

\newcommand{\bd}{\begin{deff}}

\newcommand{\ed}{\end{deff}}

\newcommand{\bl}{\begin{lm}}

\newcommand{\el}{\end{lm}}

\newcommand{\bp}{\begin{proo}}

\newcommand{\ep}{\end{proo}}

\newcommand{\bt}{\begin{thee}}

\newcommand{\et}{\end{thee}}

\newcommand{\bc}{\begin{co}}

\newcommand{\ec}{\end{co}}

\newcommand{\brm}{\begin{rem}}

\newcommand{\erm}{\end{rem}}

\newcommand{\der}{{\rm d}}

\usepackage{t1enc}
\def\frak{\mathfrak}

\newcommand{\newc}{\newcommand}

\let\ccdot\cdot
\def\cdot{\hbox to 2.5pt{\hss$\ccdot$\hss}}

\newc{\aR}{\mbox{\boldmath{$ R$}}}
\newc{\aS}{\mbox{\boldmath{$ S$}}}
\newc{\aT}{\mbox{\boldmath{$ T$}}}
\newc{\aW}{\mbox{\boldmath{$ W$}}}

\newc{\aK}{\mbox{\boldmath{$ K$}}}
\newc{\aL}{\mbox{\boldmath{$ L$}}}
\usepackage{amssymb}
\usepackage{amscd}

\let\f=\varphi

\newcommand{\bma}{\begin{pmatrix}}
\newcommand{\ema}{\end{pmatrix}}

\newcommand{\J}{{\mbox{\sf J}}}
\newc{\obstrn}[2]{B^{#1}_{#2}}

\newcommand{\rpl}                         % +) or <+
{\mbox{$
\begin{picture}(12.7,8)(-.5,-1)
\put(0,0.2){$+$}
\put(4.2,2.8){\oval(8,8)[r]}
\end{picture}$}}

\newcommand{\lpl}                         % (+ or +>
{\mbox{$
\begin{picture}(12.7,8)(-.5,-1)
\put(2,0.2){$+$}
\put(6.2,2.8){\oval(8,8)[l]}
\end{picture}$}}

\usepackage{ifthen}

\newc{\tensor}[1]{#1}
\newc{\Mvariable}[1]{\mbox{#1}}
\newc{\down}[1]{{}_{#1}}
\newc{\up}[1]{{}^{#1}}

\newc{\JulyStrut}{\rule{0mm}{6mm}}
\newc{\midtenPan}{\mbox{\sf S}}
\newc{\midten}{\mbox{\sf T}}
\newc{\midtenEi}{\mbox{\sf U}}
\newc{\ATen}{\mbox{\sf E}}
\newc{\BTen}{\mbox{\sf F}}
\newc{\CTen}{\mbox{\sf G}}

\def\sideremark#1{\ifvmode\leavevmode\fi\vadjust{\vbox to0pt{\vss% the remark
 \hbox to 0pt{\hskip\hsize\hskip1em%                          will appear only
 \vbox{\hsize3cm\tiny\raggedright\pretolerance10000%          on the side
 \noindent #1\hfill}\hss}\vbox to8pt{\vfil}\vss}}}%
\newcommand{\Span}{\mathrm{Span}}
\numberwithin{equation}{section}

\newcounter{romenumi}
\newcommand{\labelromenumi}{(\roman{romenumi})}

\newcommand{\w}{{\scriptstyle\wedge}\,}

\newcommand{\lie}{\mathcal{L}}
\newcommand{\inc}[2]{\mathbf{#1}_{\mathbf{#2}}}
\newcommand{\hc}[1]{\theta^{#1}}
\newcommand{\vc}[1]{\Omega_{#1}}

\newcommand{\bbR}{\mathbb{R}}

\newcommand{\soa}{\frak{so}}
\newcommand{\spa}{\frak{sp}}

\newcommand{\Maux}{{\text{\usefont{T1}{qcs}{m}{sl}M}}}
\newcommand{\Naux}{{\text{\usefont{T1}{qcs}{m}{sl}N}}}
\newcommand{\Waux}{{\text{\usefont{T1}{qcs}{m}{sl}W}}}

\newcommand{\vff}{\vfill\end{document}}

\begin{document}

\title[Contact projective para-CR $5$-manifolds]{Five-dimensional 
para-CR manifolds
\\
and contact projective geometry in dimension three}

\author{Jo\"el Merker} \address{}
\address{Laboratoire de Math\'ematiques d'Orsay, CNRS, 
Universit\'e Paris-Saclay, 91405 Orsay Cedex, France}
\email{joel.merker@universite-paris-saclay.fr}

\author{Pawe\l~ Nurowski} \address{Centrum Fizyki Teoretycznej,
Polska Akademia Nauk, Al. Lotnik\'ow 32/46, 02-668 Warszawa, Poland}
\email{nurowski@cft.edu.pl}

\thanks{
\\
${}\ \ \ \ \ $
This work was supported
in part by the Polish National Science Centre (NCN) 
via the grant number 2018/29/B/ST1/02583.}

\date{\today}
\begin{abstract}
We study invariant properties of $5$-dimensional para-CR structures
whose Levi form is degenerate in precisely one direction and which are
$2$-nondegenerate. We realize that {\em two}, out of three, primary
(basic) para-CR invariants of such structures are the classical
differential invariants known to {\em Monge} (1810) and
to {\em W\"unschmann} (1905):
\[
{\scriptstyle{
\Maux(G)
\,:=\,
40G_{ppp}^3-45G_{pp}G_{ppp}G_{pppp}+9G_{pp}^2G_{ppppp},
\ \ \ \ \
\Waux(H)
\,:=\,
9D^2H_r-27DH_p-18H_rDH_r+18H_pH_r+4H_r^3+54H_z.}}
\] 
The vanishing $\Maux(G) \equiv 0$ provides a local necessary and
sufficient condition for the graph of a function in the $(p,G)$-plane
to be contained in a conic, while the vanishing $\Waux(H) \equiv 0$
gives an {\sl if-and-only-if} condition for a 3\textsuperscript{rd}
order ODE to define a natural Lorentzian geometry on the space of its
solutions.

Mainly, we give a geometric interpretation of the {\em third} basic
invariant of our class of para-CR structures, the simplest one, of
lowest order, and of mixed nature $\Naux(G,H) :=
2G_{ppp}+G_{pp}H_{rr}$.  We establish that the vanishing $\Naux(G,H)
\equiv 0$ gives an {\em if-and-only-if} condition for the {\em two}
$3$-dimensional quotients of the para-CR manifold by its two canonical
integrable rank-$2$ distributions, to be equipped with contact
projective geometries.

A curious transformation between the W\"unschmann invariant and the
Monge invariant, first noted by us in a recent
publication~{\cite{Merker-Nurowski-2020}}, is also discussed, and its
mysteries are further revealed.
\end{abstract}

\maketitle

\vspace{-1.05cm}

\tableofcontents

\vspace{-1.5cm}

%%%%%%%%%%%%%%%%%%%%%%%%%%%%%%%%%%
\section{Introduction}\label{intr}
%%%%%%%%%%%%%%%%%%%%%%%%%%%%%%%%%%

The main features of the present article,
continuing our joint work~{\cite{Merker-Nurowski-2020}},
can be condensed into the following

\begin{theorem}\label{th0}
Consider a smooth $5$-dimensional para-CR structure $M^5$, 
whose Levi form is degenerate in precisely one direction, which is
$2$-nondegenerate, and which is defined as a system of two PDEs:
\[
z_y
=
G(x,y,z,z_x,z_{xx})
\quad\&\quad 
z_{xxx}
=
H(x,y,z,z_x,z_{xx}),
\quad
\text{for}\,\,
z=z(x,y),
\quad
\text{with complete integrability},
\]
in terms of two real $\mathcal{C}^\infty$ functions
%of five variables
$G = G(x,y,z,p,r)$ and $H = H(x,y,z,p,r)$ such that $G_r\equiv 0 \neq
G_{pp}$.

If one among three primary relative para-CR 
differential invariants vanishes identically:
\[
2G_{ppp}+G_{pp}H_{rr}
\equiv 
0,
\]
then the para-CR
structure defines {\em two} natural contact projective geometries
on certain {\em two} $3$-dimensional quotient spaces of $M^5$.
\end{theorem}

Concept explanations being required
to make the paper self contained, we start by briefly collecting:

\smallskip\noindent{\bf (a)}\,
basic facts about $5$-dimensional para-CR structures 
({\cite{Hill-Nurowski-2010}}; we follow exposition and notation 
from~{\cite{Merker-Nurowski-2020}});

\smallskip\noindent{\bf (b)}\,
rudiments of the theory of contact geometry of 3\textsuperscript{rd}
order ODEs ({\cite{Chern-1940}}; we follow~{\cite{Godlinski-2008,Godlinski-Nurowski-2009}}); and:

\smallskip\noindent{\bf (c)}\,
facts from the theory of contact projective structures 
({\cite{Fox-2005}}; we follow~{\cite{Godlinski-Nurowski-2009}}).

\smallskip

Then we will prove the above theorem.

%%%%%%%%%%%%%%%%%%%%%%%%%%%%%%%%%%
\section{Degenerate $5$-Dimensional Para-CR-Structures}
\label{degenerate-5-para-CR}
%%%%%%%%%%%%%%%%%%%%%%%%%%%%%%%%%%

Recall from~{\cite{Merker-2008, Hill-Nurowski-2010}} that a {\sl
para-CR structure} is a geometric structure which a hypersurface
$M^{2n-1}\subset (\bbR^n\times \bbR^n)$ acquires from the ambient
\emph{product} space $\bbR^n\times \bbR^n$. More specifically one
considers a local hypersurface
\[
M_{2n-1}
=
\big\{
\bbR^n\times\bbR^n\ni(x,\bar{x})
~|~
\Phi(x_1,\dots,x_n,\bar{x}_1,\dots,\bar{x}_n)=0
\big\},
\]
with $d_x\Phi \neq 0 \neq d_{\overline{x}} \Phi$, modulo (local)
diffeomorphisms $\varphi:\bbR^n\times\bbR^n\to \bbR^n\times\bbR^n$
preserving the splitting of $\bbR^{2n}$ into
$\bbR^{2n}=\bbR^n\times\bbR^n$,
i.e. $\varphi(x,\bar{x})=(\psi(x),\bar{\psi}(x))$, where
$\psi:\bbR^n\to\bbR^n$ and $\bar{\psi}:\bbR^n\to \bbR^n$ are (local)
diffeomorphisms.

The lowest dimension where these structures are interesting is
$n=2$. If nondegenerate, such para-CR structures are 
in 1-1 correspondence with 
2\textsuperscript{nd} order ODEs considered modulo point
transformations of variables~{\cite{
Nurowski-Sparling-2003, Merker-2008}}. In this
article we will deal with the next dimension, $n=3$, and will study
5-dimensional para-CR structures.

A 5-dimensional para-CR structure, i.e. a hypersurface $M^5\subset
\bbR^3\times\bbR^3$ considered modulo split transformations of
the product $\bbR^3\times\bbR^3$, can be defined in terms of a graph
of a function $z$ of five variables,
$z=z(x,y,\bar{x},\bar{y},\bar{z})$, where
$(x,y,z,\bar{x},\bar{y},\bar{z})$ are coordinates in
$\bbR^6=\bbR^3\times\bbR^3$. This in turn, can be considered as a
\emph{general solution} to a completely integrable
system of two PDEs on the plane $(x,y)$
for a function $z=z(x,y)$, in which $(\bar{x},\bar{y},\bar{z})$ denote
constants of integration and parametrize the solution space of the
corresponding system of PDEs.

\begin{example}
\label{Example-model}
{\bf [Model]}
Take $(x-\bar{x})^2+(y-\bar{y})(z-\bar{z})=0$, and solve it for $z$
obtaining: $z=-\frac{(x-\bar{x})^2}{y-\bar{y}}+\bar{z}$. Now think
about $(x,y)$ as independent variables, and
$(\bar{x},\bar{y},\bar{z})$ as parameters.  Obviously $z_{xxx}=0$.
Also, because $z_y=\frac{(x-\bar{x})^2}{(y-\bar{y})^2}$ and
$z_x=\frac{-2(x-\bar{x})}{(y-\bar{y})}$, we have $z_y=\tfrac14 z_x^2$.
So, a para-CR structure defined by the cone
$(x-\bar{x})^2+(y-\bar{y})(z-\bar{z})=0$ in $\bbR^3\times\bbR^3$
defines a system of PDEs on the plane
\[
\boxed{z_{xxx}=0\quad\quad\&\quad\quad z_y=\tfrac14 z_x^2\quad\quad\mathrm{for}\quad\quad z=z(x,y)}\,.
\] 

Conversely, given this system of PDEs, $z_{xxx}=0$ solves as
$z=\alpha(y)x^2+\beta(y)x+\gamma(y)$, and $z_y=\tfrac14z_x^2$ gives
sucessively: $\alpha'=\alpha^2$, hence $\alpha=\frac{-1}{y-\bar{y}}$,
$\beta'=\frac{-\beta}{y-\bar{y}}$, hence
$\beta=\frac{2\bar{x}}{y-\bar{y}}$,
$\gamma'=\frac{\bar{x}{}^2}{(y-\bar{y})^2}$, hence
$\gamma=\frac{-\bar{x}{}^2}{y-\bar{y}}+\bar{z}$.  This finally gives
$z=\frac{-\bar{x}{}^2}{y-\bar{y}}+\bar{z}+\frac{2x\bar{x}}{y-\bar{y}}-\frac{x^2}{y-\bar{y}}$,
i.e. the cone
\[
\boxed{(x-\bar{x})^2+(y-\bar{y})(z-\bar{z})=0}\,.
\eqno\qed
\]
\end{example}

\smallskip

In general, we consider the following system of two PDEs on the plane
\begin{equation}
\boxed{
z_{xxx}
=
H(x,y,z,z_x,z_{xx})
\quad\&\quad 
z_y
=
G(x,y,z,z_x,z_{xx})
\quad\mathrm{for}
\quad z=z(x,y)}\,.\label{ss}
\end{equation}

\begin{lemma}
{\rm {\cite{Merker-2008}}}
The general solution of~{\eqref{ss}} depends on 3 parameters
$(\bar{x}, \bar{y}, \bar{z})$, and has the form $z =
z(x,y;\bar{x},\bar{y},\bar{z})$ if and only if
\begin{equation}
\boxed{\triangle H=D^3G}\,,\label{ic}
\end{equation}
where, abbreviating $p=z_x$, $r=z_{xx}$,
\[
D
=
\partial_x+p\partial_z+r\partial_p+H\partial_r,
\quad\quad
\triangle
=
\partial_y+G\partial_z+DG\partial_p+D^2G\partial_r.
\eqno\qed
\] 
\end{lemma}

General solutions of systems~{\eqref{ss}} give examples of 5-dimensional
para-CR structures.  We prefer the PDE point of view, and we will
stick to this in the following.  In particular, in this point of view,
para-CR transformations for hypersurfaces in
$(x,y,z,\bar{x},\bar{y},\bar{z})$, are the \emph{point transformations
of variables} of~{\eqref{ss}}.

Thus, we can either describe our para-CR geometry as a geometry of
hypersurfaces in the $(x,y,z,\bar{x},\bar{y},\bar{z})$ space (modulo
appropriate diffeomorphisms), or as a geometry of PDEs~{\eqref{ss}}
considered modulo point transformation of variables.

It is clear from the hypersurfaces picture, that a 5-dimensional
para-CR manifold $M^5$ is equipped with \emph{two integrable
distributions} $D_1$ and $D_2$. These are tangent to the
\emph{foliations} of $M^5$ obtained by intersecting it with either the
$3$-planes $\{x=\mathrm{const}, y=\mathrm{const}, z=\mathrm{const}\}$,
or the $3$-planes $\{\bar{x}=\mathrm{const}, \bar{y}=\mathrm{const},
\bar{z}=\mathrm{const}\}$.

In the PDE picture, these two distributions are the respective
\emph{annihilators} of the following system of 1-forms
\begin{equation}
\begin{array}{lll}
D_1=\bma
\omega^1=\der z-p\der x-G\der y\\ \omega^2=\der p-r \der x-DG\der y\\\omega^3=\der r-H\der x-D^2G\der y \ema^\perp &\&&
D_2=\bma
\omega^1=\der z-p\der x-G\der y\\ \omega^4=\der x\\\omega^5=\der y 
\ema^\perp.
\end{array}
\label{d1d2}
\end{equation}

Actually, the condition that $D_1$ is integrable is precisely the integrability condition~{\eqref{ic}} guaranteeing that the PDE system~{\eqref{ss}} has a 3-parameter family of solutions~{\cite{Merker-2008}}.
Note that the rank $4$ distribution $D=D_1+D_2$ 
is also well defined.

This enables for a {\em definition} of a $5$-dimensional para-CR
structure, locally, `{\sl \`a la \'Elie Cartan}'.

\begin{definition}
{\em A 5-dimensional para-CR structure is a structure consisting of
an equivalence class $[\omega]$ of coframes
$\omega=(\omega^1,\omega^2,\omega^3,\omega^4,\omega^5)$ on $\bbR^5$
parameterized by $(x,y,z,p,r)$, with an equivalence relation $\sim$
given by}
\[
\bar{\omega}\sim\omega
\quad\quad\Longleftrightarrow\quad\quad 
\bma
\bar{\omega}^1\\\bar{\omega}^2\\\bar{\omega}^3\\\bar{\omega}^4\\\bar{\omega}^5\ema=\bma f_1&0&0&0&0\\f_2&\rho\mathrm{e}^\phi&f_4&0&0\\f_5&f_6&f_7&0&0\\\bar{f}_2&0&0&\rho\mathrm{e}^{-\phi}&\bar{f}_4\\\bar{f}_5&0&0&\bar{f}_6&\bar{f}_7
\ema
\bma\omega^1\\\omega^2\\\omega^3\\\omega^4\\\omega^5
\ema,
\]
{\em with $\omega^1=\der z-p\der x-G\der y$, $\omega^2=\der p-r \der
x-DG\der y$, $\omega^3=\der r-H\der x-D^2G\der y$, $\omega^4=\der x$,
$\omega^5=\der y$, being in the class $[\omega]$.}
\end{definition}

The integrabilities of the two distributions $D_1$ and $D_2$, 
as defined in~{\eqref{d1d2}}, implies that
\[
\Big(\der \omega^1
-L_{11}\omega^2\w\omega^4
-L_{12}\omega^2\w\omega^5
-L_{21}\omega^3\w\omega^4
-L_{22}\omega^3\w\omega^5\Big)
\w
\omega^1\equiv 0,
\]
with a certain $2\times 2$ matrix $L$ of functions $L_{AB}$,
$A,B=1,2$, on $M^5$ defined by this condition.

The matrix $L$, called the \emph{Levi form}, is \emph{not} well
defined by the equivalence class of $\omega$, but \emph{its signature
is}. Hence $\det(L)=0$, or $\det(L)\neq 0$, is a para-CR invariant
condition at each point. If $\det(L)\neq 0$, the corresponding para-CR
structure is \emph{non}degenerate, and it defines one of the
\emph{parabolic} geometries in dimension 5 (flat model\,\,---\,\,a
flying soucer in the \emph{attacking mode}).

\medskip

In this paper, we consider para-CR structures with
\[
L\neq 0\quad\quad 
\mathrm{but\,\,such\,\,that}
\quad\quad\det(L)\equiv 0.
\]
These are 5-dimensional para-CR structures 
\emph{with Levi form $L$ degenerate in 1 direction}.

\medskip

In terms of our PDEs, this degeneracy means that
\begin{equation}
G_r\equiv 0,
\quad\text{that is}\quad 
G=G(x,y,z,z_x).
\label{1deg}
\end{equation}
We also \emph{do not want} that our para-CR structure is locally
para-CR-equivalent to a product of a 3-dimensional para-CR manifold
$M^3$ and a product $\bbR\times \bbR$. This results in our further
assumption that 
\begin{equation}
G_{pp}\neq 0.
\label{2ndeg}
\end{equation}

%%%%%%%%%%%%%%%%%%%%%%%%%%%%%%%%%%%%%%%%%%%%%%%%%%%%%%%%%%%%
\section{Basic invariants for Degenerate para-CR Structures}
%%%%%%%%%%%%%%%%%%%%%%%%%%%%%%%%%%%%%%%%%%%%%%%%%%%%%%%%%%%%

Summarizing, we study systems of PDEs on the plane:
\[
\boxed{z_{xxx}=H(x,y,z,p,r)
\quad\quad\&\quad\quad 
z_y=G(x,y,z,p)
\quad\quad\mathrm{for}\quad\quad z(x,y)}\,,
\]
such that
\[
\boxed{\triangle H=D^3G}\quad\quad\&\quad\quad\boxed{G_{pp}\neq 0}\,,
\]
with $D=\partial_x+p\partial_z+r\partial_p+H\partial_r$,
$\triangle=\partial_y+G\partial_z+DG\partial_p+D^2G\partial_r$, and
$p=z_x$, $r=z_{xx}$, considered modulo \emph{point transformations of
variables}. This is equivalent to study coframes $\omega^1=\der
z-p\der x-G\der y$, $\omega^2=\der p-r \der x-DG\der y$,
$\omega^3=\der r-H\der x-D^2G\der y$, $\omega^4=\der x$,
$\omega^5=\der y$, with $D^3G=\triangle H$, $G_{pp}\neq 0$, and
$G_r \equiv 0$, given modulo
\begin{equation}
\bma
\omega^1\\\omega^2\\\omega^3\\\omega^4\\\omega^5
\ema
\longmapsto
\bma 
f_1&0&0&0&0\\f_2&\rho\mathrm{e}^\phi&f_4&0&0\\f_5&f_6&f_7&0&0\\\bar{f}_2&0&0&\rho\mathrm{e}^{-\phi}&\bar{f}_4\\\bar{f}_5&0&0&\bar{f}_6&\bar{f}_7\ema\bma\omega^1\\\omega^2\\\omega^3\\\omega^4\\\omega^5
\ema.
\label{trump}
\end{equation}
In reference~{\cite{Merker-Nurowski-2020}}, studying such structures,
we established among other things, the following

\begin{theorem}\label{the1}
It is always possible to invariantly force the lifted coframe
$\theta^1=f_1\omega^1$,
$\theta^2=f_2\omega^1+\rho\mathrm{e}^\phi\omega^2+f_4\omega^3$,
$\theta^3=f_5\omega^1+f_6\omega^2+f_7\omega^3$,
$\theta^4=\bar{f}_2\omega^1+\rho\mathrm{e}^{-\phi}\omega^4+\bar{f}_4\omega^5$,
$\theta^5=\bar{f}_5\omega^1+\bar{f}_6\omega^2+\bar{f}_7\omega^3$ to
satisfy the following EDS:
\[
\begin{aligned}
\der\hc{1} 
=&~
\vc{1}\w\hc{1}+\hc{2}\w\hc{4},\nonumber\\
\der\hc{2} 
=&~
\hc{2}\w(\vc{2}-\tfrac12\vc{1})-\hc{1}\w\vc{3}+\hc{3}\w\hc{4},\nonumber \\
\der\hc{3}
=&~
2\hc{3}\w\vc{2}-\hc{2}\w\vc{3}+Q\hc{1}\w\hc{3}-\tfrac12(\tfrac{\mathrm{e}^{\phi}}{3\rho})^3{\color{red}\inc{A}{}}\hc{1}\w\hc{4}+\tfrac{\mathrm{e}^{-\phi}}{3\rho}\inc{{\color{blue}C}}{}\hc{2}\w\hc{3}, \nonumber\\
\der\hc{4}
=&~
-\hc{2}\w\hc{5}-\hc{4}\w(\tfrac12\vc{1}+\vc{2})-\hc{1}\w\vc{4},\nonumber \\
\der\hc{5}
=&~
-2\hc{5}\w\vc{2}+\hc{4}\w\vc{2}+(\tfrac{\mathrm{e}^{\phi}}{3\rho})^3{\color{green}\inc{B}{}}\hc{1}\w\hc{2}+Q\hc{1}\w\hc{5}+\tfrac{\mathrm{e}^{\phi}}{3\rho}\tilde{\inc{C}{}}\hc{4}\w\hc{5},
\end{aligned}
\]
in which three primary relative differential invariants are
\[
\begin{aligned}
{\color{red}\inc{A}{}}
&=
{\color{red}9D^2H_r-27DH_p-18H_rDH_r+18H_pH_r+4H_r^3+54H_z},\\
{\color{green}\inc{B}{}}
&=
(\tfrac{1}{2G_{pp}^3})~[~{\color{green}40G_{ppp}^3-45G_{pp}G_{ppp}G_{pppp}+9G_{pp}^2G_{ppppp}}~],\\
\inc{{\color{blue}C}}{}
&=
(\tfrac{1}{G_{pp}})~[~{\color{blue}2G_{ppp}+G_{pp}H_{rr}}~],
\end{aligned}
\]
that is, the vanishing or not of each of ${\color{red}\inc{A}{}}$, 
${\color{green} \inc{B}{}}$, $\inc{{\color{blue}C}}{}$ 
is \emph{an invariant property} 
of the corresponding para-CR structure.
Lastly, $\tilde{ \mathbf{C}}$ vanishes identically when 
$\inc{{\color{blue} C}}{} \equiv 0$.
\end{theorem}

\begin{remarks}$\:$

\begin{itemize}

\item 
\emph{Flat model}:
${\color{red}\inc{A}{}} = {\color{green}\inc{B}{}} =
\inc{{\color{blue}C}}{}=0$, and this is locally equivalent to
$z_{xxx}=0$, $z_y=\tfrac14z_x^2$, i.e. to the para-CR structure from
our Example~{\ref{Example-model}} in the beginning,
{\em cf.}~{\cite{Merker-Nurowski-2020}}.

\smallskip\item 
\emph{Symmetries}: A vector field $X$ on $M^5\ni(x,y,z,p,r)$ is a \emph{symmetry} of the para-CR structure as defined 
in~{\eqref{ss}}--{\eqref{d1d2}} if and only if
\[
\begin{aligned}
&\big(\lie_X\omega^1\big)\w\omega^1=0,\\
&\big(\lie_X\omega^2\big)\w\omega^1\w\omega^2\w\omega^3=0,\quad\quad  &\big(\lie_X\omega^4\big)\w\omega^1\w\omega^4\w\omega^5=0,\\
&\big(\lie_X\omega^3\big)\w\omega^1\w\omega^2\w\omega^3=0, \quad\quad  &\big(\lie_X\omega^5\big)\w\omega^1\w\omega^4\w\omega^5=0.
\end{aligned}
\]
Any Lie bracket of two symmetries is a symmetry, which brings the
notion of a \emph{symmetry algebra} of a para-CR-structure: the Lie
algebra over the reals of all symmetries.

\smallskip\item 
For our flat model with ${\color{red} \inc{A}{}} = {\color{green}
\inc{B}{}}=\inc{{\color{blue} C}}{}=0$, the \emph{symmetry algebra}
is $\spa(4,\bbR)\simeq \soa(2,3)$.
\end{itemize}

\end{remarks}

%%%%%%%%%%%%%%%%%%%%%%%%%%%%%%%%%%%%%%%%%%%%%%%%%%%%%%%
\section{Geometry of W\"unschmann and Monge Invariants}
%%%%%%%%%%%%%%%%%%%%%%%%%%%%%%%%%%%%%%%%%%%%%%%%%%%%%%%

The explicit expressions for the \emph{relative invariants}
${\color{red}\inc{A}{}}$ and ${\color{green}\inc{B}{}}$ of the
considered para-CR structures redirect us to the
\emph{{\color{red}theory of 3\textsuperscript{rd} order ODEs
considered modulo contact transformations of variables}} and to
\emph{{\color{green}differential geometry of conics on the plane}}. We
therefore make the following interlude in our main theme now.

\subsection{3\textsuperscript{rd} order ODEs 
considered modulo contact transformation of variables}
\label{3-rd-order-ODE-point-transformations} 
We formulate a theorem~{\cite{Godlinski-2008,
Godlinski-Nurowski-2009}} about the main structure which is
associated with third-order ODEs modulo contact transformations of
variables, namely about an $\spa(4,\bbR)$-valued Cartan connection on
the bundle $P^{10}\to\J^2$. This structure will serve as a starting
point for analyzing further geometries of ODEs.

\begin{theorem}\label{th.c.1}
To every third order ODE $z'''=H(x,z,z',z'')$, there is associated a
(principal) fibre bundle $H_6\to P^{10}\to\J^2$, over the space of
second jets, where $\dim P^{10}=10$ and $H_6$ is an appropriate
six-dimensional subgroup of $SP(4,\bbR)$, with the group parameters
$u_i$, $i=1,2\dots, 6$, and a unique coframe of 1-forms $(\theta^1$,
$\theta^2$, $\theta^3$, $\theta^4,\theta^5$, $\Omega_1$, $\Omega_2$,
$\Omega_3$, $\Omega_4$, $\Omega_5)$ on $P^{10}$, which satisfies the
following EDS:
\begin{align}
\der\hc{1}
=&~
\vc{1}\w\hc{1}+\hc{4}\w\hc{2},
\nonumber \\
\der\hc{2}
=&~
\vc{2}\w\hc{1}+\vc{3}\w\hc{2}+\hc{4}\w\hc{3},
\nonumber \\
\der\hc{3}
=&~
\vc{2}\w\hc{2}+(2\vc{3}-\vc{1})\w\hc{3}+\inc{A}{2}\hc{2}\w\hc{1}+{\color{red}\inc{A}{1}}\hc{4}\w\hc{1}, 
\nonumber \\
\der\hc{4}
=&~
\vc{4}\w\hc{1}+(\vc{1}-\vc{3})\w\hc{4}+\hc{5}\w\hc{2},
\nonumber \\
\der\hc{5}
=&~
\vc{4}\w\hc{4}
+(\vc{1}-2\vc{3})\w\hc{5}+(\inc{A}{7}+\inc{Z}{3})\hc{1}\w\hc{2}+\inc{Z}{4}\hc{1}\w\hc{3} 
\nonumber \\
&-\inc{A}{5}\hc{1}\w\hc{4}+{\color{blue}\inc{Z}{1}}\hc{2}\w\hc{3},
\nonumber \\
\der\vc{1}
=&~
\vc{5}\w\hc{1}+\vc{4}\w\hc{2}-\vc{2}\w\hc{4},
\nonumber \\
\der\vc{2}
=&~
(\vc{3}-\vc{1})\w\vc{2}+\tfrac{1}{2}\vc{5}\w\hc{2}+\vc{4}\w\hc{3}
+\inc{A}{3}\,\hc{1}\w\hc{2}+\inc{A}{4}\hc{1}\w\hc{4}, 
\nonumber \\
\der\vc{3}
=&~
\tfrac{1}{2}\vc{5}\w\hc{1}+\vc{4}\w\hc{2}+\hc{5}\w\hc{3}
+\inc{A}{5}\hc{1}\w\hc{2}+\inc{A}{2}\hc{1}\w\hc{4},
\label{e.c.dtheta_10d} \\
\der\vc{4}
=&~
\hc{5}\w\vc{2}+\vc{4}\w\vc{3}+\tfrac{1}{2}\vc{5}\w\hc{4}
+(\inc{A}{6}+\inc{Z}{2})\hc{1}\w\hc{2} +2\inc{Z}{3}\hc{1}\w\hc{3} 
\nonumber, \\
&-\inc{A}{3}\hc{1}\w\hc{4}+\inc{Z}{4}\hc{2}\w\hc{3} 
\nonumber \\
\der\vc{5}
=&~
\vc{5}\w\vc{1}+2\vc{4}\w\vc{2}+\inc{C}{1}\hc{1}\w\hc{2}
+2\inc{Z}{2}\hc{1}\w\hc{3}+\inc{A}{8}\hc{1}\w\hc{4}+2\inc{Z}{3}\hc{2}\w\hc{3}.
\nonumber
\end{align}
Here ${\color{red} \inc{A}{1}}, \ldots, \inc{A}{8}, {\color{blue}
\inc{Z}{1}}, \ldots, \inc{Z}{4}, \inc{C}{1}$ are functions on
$P^{10}$.

The $8 + 4 + 1$ functions ${\color{red} \inc{A}{1}}, \ldots,
{\color{blue}\inc{Z}{1}}, \dots, \inc{C}{1}$ are contact relative
invariants of the underlying ODE and the full set of contact
invariants can be constructed by consecutive differentiations of
${\color{red} \inc{A}{1}}, \ldots, {\color{blue} \inc{Z}{1}}, \dots,
\inc{C}{1}$ with respect to the frame $(X_1,X_2,X_3,X_4,X_5,X_6,\,
X_7,X_8,X_9, X_{10})$ dual to $(\theta^1, \theta^2, \theta^3,
\theta^4, \Omega_1, \Omega_2, \Omega_3, \Omega_4, \Omega_5,
\Omega_6)$.

The coframe $(\theta^1, \theta^2, \theta^3, \theta^4, \Omega_1,
\Omega_2, \Omega_3, \Omega_4, \Omega_5, \Omega_6)$ defines the
$\spa(4, \bbR)$-valued Cartan normal connection $\wh{\omega}$ on
$P^{10}$ by
\begin{equation}\label{e.c.conn_sp}
\wh{\omega}
=
\bma
\tfrac{1}{2}\vc{1} & \tfrac{1}{2}\vc{2} & -\tfrac{1}{2}\vc{4} & -\tfrac{1}{4}\vc{5} 
\\\\
\hc{4} & \vc{3}-\tfrac{1}{2}\vc{1} & -\hc{5} & -\tfrac{1}{2}\vc{4} 
\\\\
\hc{2} & \hc{3} & \tfrac{1}{2}\vc{1}-\vc{3} & -\tfrac{1}{2}\vc{2} 
\\\\
2\hc{1} & \hc{2} & -\hc{4} & -\tfrac{1}{2}\vc{1}
\ema.
\end{equation}
The EDS~{\eqref{e.c.dtheta_10d}} gives explicit formulas for the
\emph{curvature} $\wh{K} = \der \wh{\omega} +
\wh{\omega}\w\wh{\omega}$ of this Cartan normal connection, with the
invariant functions $\inc{A}{\alpha}, \inc{Z}{\beta},\inc{C}{1}$,
being the apropriate entries in the coframe components matrices
$\wh{K}_{ij}$ of $\wh{K} = \tfrac12 \wh{K}_{ij}\theta^i\w\theta^j$.

Two 3\textsuperscript{rd} order ODEs $y'''=F(x,y,y',y'')$ and
$\bar{y}'''=\bar{F}(\bar{x}, \bar{y}, \bar{y}',
\bar{y}'')$ are locally contact equivalent if and
only if their associated Cartan connections are locally
diffeomorphic, that is, there exists a local bundle diffeomorphism
$\Phi\colon\bar{P}\to P$ such that
\[
\Phi^*\wh{\omega}
=
\overline{\wh{\omega}}.
\]
\end{theorem}

It further follows that: 

\smallskip\noindent$\bullet$\,
$\inc{A}{2},\ldots,\inc{A}{8}$ express in terms of 
coframe derivatives of ${\color{red}\inc{A}{1}}$; 

\smallskip\noindent$\bullet$\,
$\inc{Z}{2},\ldots,\inc{Z}{4}$ express 
in terms of coframe derivatives of
${\color{blue}\inc{Z}{1}}$;

\smallskip\noindent$\bullet$\,
$\inc{C}{1}$ is a function of
coframe derivatives of both ${\color{red}\inc{A}{1}}$ and
${\color{blue}\inc{Z}{1}}$. 

\smallskip\noindent
So only ${\color{red}\inc{A}{1}}$ and
${\color{blue}\inc{Z}{1}}$ are \emph{basic} 
(primary) invariants, namely all
other (secondary) invariants are deduced by differentiation. 
Their remarkable explicit expressions are given by

\begin{proposition}\label{c.cor.har}
Letting $D=\partial_x+p\partial_z+r\partial_p+H\partial_r$, 
and $u_1$ and $u_3$ be the parameters along 
the gauge group $H_6$ mentioned in 
Theorem~{\ref{th.c.1}}, one has:
\[
\begin{aligned}
{\color{red}\inc{A}{1}}
=&~
\tfrac12(\frac{u_3}{3u_1})^3~[~{\color{red}9D^2H_r-27DH_p-18H_rDH_r+18H_pH_r+4H_r^3+54H_z}~]
&=:&~
\tfrac12(\frac{u_3}{3u_1})^3~{\color{red}\inc{A}{}},
\\
{\color{blue}\inc{Z}{1}}
=&~
\frac{u_1^2}{6u_3^5}~{\color{blue}H_{rrrr}}
&=:&~
\frac{u_1^2}{6u_3^5}~{\color{blue}\inc{Z}{}}.
\end{aligned}
\]
\end{proposition}

Thus, the \emph{contact} relative invariant ${\color{red}\inc{A}{1}}$
for a contact equivalence class of ODEs $z'''=H(x,z,z',z'')$ is given,
modulo a nonvanishing scaling factor, by the same expression as one of
our basic para-CR invariants ${\color{red}\inc{A}{}}$ for the
5-dimensional para-CR manifolds with Levi form degenerate in one
direction\footnote{It is not a big surprise, though, since our PDEs on
the plane~{\eqref{ss}} include a one parameter family of ODEs
$z'''=H(x,y,z,z',z'')$, parametrized by the variable $y$.}.

\medskip

The expression
${\color{red}\inc{A}{}} = 
{\color{red}9D^2H_r-27DH_p-18H_rDH_r+18H_pH_r+4H_r^3+54H_z}$
was for the first time obtained in 1905 by
W\"unschmann~{\cite{Wunschmann-1905}}, who observed that its vanishing
or not is a contact invariant property of an ODE
$z''' = H(x,z,z',z'')$. More importantly, he also established the
\emph{geometric interpretation} of the \emph{vanishing} of
${\color{red} \inc{A}{}}$. According to W\"unschmann, if
${\color{red} \inc{A}{} \equiv 0}$, the 3-dimensional solution space of
the ODE $z'''=H(x,z,z',z'')$ is naturally equipped with a conformal
Lorentzian structure; moreover, there is a \emph{local one-to-one
correspondence} between \emph{3-dimensional conformal Lorentzian
structures} and \emph{contact equivalence classes of ODEs
$z'''=H(x,z,z',z'')$ satisfying ${\color{red}\inc{A}{}\equiv 0}$}.

\medskip

The first person who observed that the vanishing or not of
${\color{blue}\inc{Z}{} = H_{rrrr}}$ is a contact invariant property of
the ODE $z''' = H(x,z,z',z'')$ was Chern in 1940~{\cite{Chern-1940}}.
The \emph{geometric meaning of the condition that
${\color{blue}\inc{Z}{}}$ vanishes} is less
known~{\cite{Godlinski-2008}}. To fully apreciate it, one needs a
rather recent notion of a \emph{contact projective
structure}~{\cite{Fox-2005}}. Here is its definition, adapted to our
case of a 3-dimensional manifold of first jets $J^1$ of the equation
$z'''=H(x,z,z',z'')$.

\begin{definition} 
A {\sl contact projective structure}
on the first jet space $J^1 \ni (x,z,p)$ consists of:
\begin{itemize}

\smallskip\item[i)] 
the contact distribution $\C$, that is the distribution annihilated by
$\omega^1=\der z -p\der x$; together with:

\smallskip\item[ii)] 
a family of unparameterized curves in $J^1$, which are everywhere
tangent to $\C$ and such that: 

\begin{itemize}

\smallskip\item[a)] 
for any given point and direction in $\C$, there is exactly one curve
passing through that point and tangent to that direction;

\smallskip\item[b)] 
curves of the family are among unparameterized
geodesics for some linear connection on $J^1$.
\end{itemize}

\end{itemize}
\end{definition}

In other words, the idea of this geometry in the context of ODEs is as
follows\footnote{Here we quote from the PhD
Thesis~{\cite{Godlinski-2008}} of Godli\'nski, who was the first to
observe this.}: \emph{Consider the solutions of the ODE
$z'''=H(x,z,z',z'')$ as a family of curves in $J^1$ and ask whether
these curves are among geodesics of a linear connection. The answer
to this question is positive if and only if
${\color{blue}H_{rrrr}\equiv 0}$, and in this case there is a whole
family of connections for which the solutions are geodesics.}

This information about the W\"unschmann, ${\color{red}\inc{A}{}}$, and
the Chern, ${\color{blue}\inc{Z}{}}$, invariants can be nicely phrased
in terms of the natural \emph{double fibration}
\begin{align}
\xymatrix{
&\mathrm{J^2} \ar[dl]_{\pi_2} \ar[dr]^{\pi_1} & \\
\mathrm{S} & & \mathrm{J^1} }
\end{align}
of the space of second jets for the ODE $z'''=H(x,z,z',z'')$ over (a)
the solution space $S$ and (b) the space of first jets $J^1$. Here,
$\pi_1$ is the natural projection from $J^2$ to $J^1$,
$\pi_1(x,z,z',z'')=(x,z,z')$, and $\pi_2$ is a projection from $J^2$
to the space of solutions $S$ identifying points on the integral
curves of the total differential vector field
$D=\partial_x+z'\partial_z+z''\partial_{z'}+H\partial_{z''}$ on
$J^2$. In terms of this double fibration, we have the following
proposition,
in which $z'=p$,
$z''=r$, and $D=\partial_x+p\partial_z+r\partial_p+H\partial_r$.

\begin{proposition}\label{33}
Two basic (primary) local contact relative invariants for
3\textsuperscript{rd} order ODEs $z'''=H(x,z,z',z'')$
are the \emph{W\"unschmann invariant}, ${\color{red}\inc{A}{}}
=~{\color{red}9D^2H_r-27DH_p-18H_rDH_r+18H_pH_r+4H_r^3+54H_z}$, and
the Chern invariant, ${\color{blue}\inc{Z}{}=H_{rrrr}}$.

The vanishing of the W\"unschmann invariant,
${\color{red}\inc{A}{}\equiv 0}$, is equivalent to have a
\emph{conformal Lorentzian structure} on the solution space $S$,
while the vanishing of the Chern invariant,
${\color{blue}\inc{Z}{}\equiv 0}$, is equivalent to have a
\emph{contact projective structure} on the space of first jets
$J^1$.
\end{proposition}

%%%%%%%%%%%%%%%%%%%%%%%%%%%%%%%%
\subsection{Conics on the plane}
Consider the most general \emph{conic} on the plane $\bbR^2$
parameterized by $(p,G)\in\bbR^2$. Such a conic is a curve in $\bbR^2$
given by the equation
\[
a_1G^2+2a_2 p G+a_3 p^2+a_4 G+a_5 p+a_6
=
0,
\]
and $a_1,\ldots, a_6$ are real constants. One can think about the
equation $a_1G^2+2a_2 p G+a_3 p^2+a_4 G+a_5 p+a_6=0$ as an implicit
relation for a function $G=G(p)$, whose graph on the plane is a conic.
It was Monge~{\cite{Monge-1810}}, who in 1810 found a differential
equation satisfied by this function. To get this equation one
eliminates $a_2,\ldots,a_6$ from the system of linear equations
\[
\frac{\der^k}{\der p^k}
\Big(a_1G^2+2a_2 p G+a_3 p^2+a_4 G+a_5 p+a_6\Big)=0, \quad\mathrm{for\,\,all}\quad k=0,1,2,3,4,5.
\]
The result is
\[
a_1G_{pp}~(~{\color{green}
40 G_{ppp}^3-45G_{pp}G_{ppp}G_{pppp}+9G_{pp}^2G_{ppppp}}~)
=
0.
\]
Excluding the nongeneric case when $a_1G_{pp}=0$, one obtains the Monge 5th order ODE
\[
{\color{green}
40 G_{ppp}^3-45G_{pp}G_{ppp}G_{pppp}+9G_{pp}^2G_{ppppp}}
=
0
\]
for a local function $G=G(p)$ to have a graph contained
in a general conic. 

In the context of this paper it is necessary to note, that the left
hand side of this expression
${\color{green}\inc{M}{}}:={\color{green}40
G_{ppp}^3-45G_{pp}G_{ppp}G_{pppp}+9G_{pp}^2G_{ppppp}}$ is, modulo a
nonvanishing factor, the same as the relative para-CR invariant
${\color{green}\inc{B}{}}$ for 5-dimensional para-CR structures given
by~{\eqref{ss}}--{\eqref{ic}}, {\eqref{1deg}}--{\eqref{2ndeg}}. More
precisely, the vanishing of ${\color{green}\inc{B}{}}$ is equivalent
to the vanishing of a 3-parameter family of Monge 5th order ODEs
${\color{green}\inc{M}{}=0}$, with parameters $x,y,z$.

This justifies our terminology, which we adopt from now on, that the
relative para-CR invariant 
\[
{\color{green}\inc{B}{}}=\frac{1}{2G_{pp}^3}{\color{green}\inc{M}{}},
\]
or its core
\[
{\color{green}\inc{M}{}}={\color{green}40
G_{ppp}^3-45G_{pp}G_{ppp}G_{pppp}+9G_{pp}^2G_{ppppp}},
\]
will be called the \emph{Monge invariant}.

In this way we have a nice geometric interpretation of the vanishing
of the para-CR invariant ${\color{green}\inc{B}{}}$: it vanishes if
and only if $G=G(x,y,p,z)$ defines a (general) 
conic on the plane $(p,G)$.

We close this section with a remark that we have yet another geometric
interpretation of the vanishing of the invariant
${\color{green}\inc{B}{}}$. This is described in our recent
paper~{\cite{Merker-Nurowski-2019}}, and is related to the single PDE
$z_y=G(x,y,z,z_x)$ for a function $z=z(x,y)$ considered modulo point
transformations of variables.

%%%%%%%%%%%%%%%%%%%%%%%%%%%%%%%%%%%%%%%%%
\section{5-Dimensional Para-CR Structures 
as 3\textsuperscript{rd} Order ODEs}
\label{5-para-CR-3rd-ODE}
%%%%%%%%%%%%%%%%%%%%%%%%%

Theorem~{\ref{the1}}, which we invoked in Section 2 of the present
paper, has its more technical, but also more refined, version which we
need now. We quote it from reference~{\cite{Merker-Nurowski-2020}}.

%/home/pawel/niebieski/notebooks/merker/homogeneous models/da_capo_2019/december_2019/start_abstract_dalej.nb
%.../start_abstract_dalej_2_Cartan_connection_Ineq0.nb
\begin{theorem}\label{susend}
Given the 1-forms
\[
\begin{aligned}
\omega^1=&\der z-p\der x-G\der y,\\
\omega^2=&\der p-r\der x-DG\der y,\\\omega^3=&\der r -H\der x-D^2 G\der y,\\
\omega^4=&\der x,\quad\quad\omega^5=\der y,
\end{aligned}
\]
 representing a $5$-dimensional para-CR manifold with $G_r=0$ and $G_{pp}\neq 0$ one can always find a para-CR equivalent set of 1-forms
\[
\begin{aligned}
\bar{\om}^1
=&~
f_1\om^1,\\
\bar{\om}^2
=&~
f_2\om^1+\rho\mathrm{e}^\phi\om^2+f_4\om^3,\\
\bar{\om}^3
=&~
f_5\om^1+f_6 \om^2+f_7\om^3,\\
\bar{\om}^4
=&~
\bar{f}_2\om^1+\rho\mathrm{e}^{-\phi}\om^4+\bar{f}_4\om^5,\\
\bar{\om}^5
=&~
\bar{f}_5\om^1+\bar{f}_6\om^4+\bar{f}_7\om^5,
\end{aligned}
\]
and additional 1-forms $\varpi_1$, $\varpi_2$, $\varpi_3$, $\varpi_4$
with
\[
\bar{\om}^1\w \bar{\om}^2\w \bar{\om}^3\w \bar{\om}^4\w \bar{\om}^5\w \varpi_1\w \varpi_2\w\varpi_3\w\varpi_4\neq 0,
\]
such that the nine 1-forms $(\bar{\om}^1,\bar{\om}^2,\bar{\om}^3,
\bar{\om}^4,\bar{\om}^5, \varpi_1, \varpi_2,\varpi_3,\varpi_4)$
satisfy the following EDS:
\begin{equation}
\begin{aligned}
\der \bar{\om}^1
=&~
-\bar{\om}^1\dz\varpi_1+\bar{\om}^2\dz\bar{\om}^4,\\
\der \bar{\om}^2
=&~
-\bar{\om}^1\dz\varpi_3+\bar{\om}^2\dz(\varpi_2-\tfrac12\varpi_1)+\bar{\om}^3\dz\bar{\om}^4,\\
\der \bar{\om}^3
=&~
-\bar{\om}^2\dz\varpi_3+2\bar{\om}^3\dz\varpi_2+\tfrac18(2I^3{}_{|4}+I^3{}_{|52})\bar{\om}^1\dz\bar{\om}^3+\\
&\ \ \ 
{\color{red}I^1}\,\bar{\om}^1\dz\bar{\om}^4+{\color{blue}I^3}\,\bar{\om}^2\dz\bar{\om}^3,\\
\der \bar{\om}^4
=&~
-\bar{\om}^1\dz\varpi_4-\bar{\om}^4\dz(\varpi_2+\tfrac12\varpi_1)-\bar{\om}^2\dz\bar{\om}^5,\\
\der \bar{\om}^5
=&~
\bar{\om}^4\dz\varpi_4-2\bar{\om}^5\dz\varpi_2+{\color{green}I^2}\,\bar{\om}^1\dz\bar{\om}^2+\tfrac18(2I^3{}_{|4}+I^3{}_{|52})\bar{\om}^1\dz\bar{\om}^5-\\
&~\tfrac12I^3{}_{|5}\,\bar{\om}^4\dz\bar{\om}^5.
\end{aligned}\label{sysend}
\end{equation}
\[
\begin{aligned}
\der {\color{red}I^1}
=&\,
I^1{}_{|1}\bar{\om}^1+I^1{}_{|2}\bar{\om}^2+I^1{}_{|3}\bar{\om}^3+I^1{}_{|4}\bar{\om}^4-\tfrac32 {\color{red}I^1}\varpi_1-3{\color{red}I^1}\varpi_2,\\
\der {\color{green}I^2}
=&\,
I^2{}_{|1}\bar{\om}^1+I^2{}_{|2}\bar{\om}^2+I^2{}_{|4}\bar{\om}^4+I^2{}_{|5}\bar{\om}^5-\tfrac32 {\color{green}I^2}\varpi_1+3{\color{green}I^2}\varpi_2,\\
\der {\color{blue}I^3}
=&\,
I^3{}_{|1}\bar{\om}^1+I^3{}_{|2}\bar{\om}^2+I^3{}_{|3}\bar{\om}^3+I^3{}_{|4}\bar{\om}^4+I^3{}_{|5}\bar{\om}^5-\tfrac12 {\color{blue}I^3}\varpi_1+ {\color{blue}I^3}\varpi_2,
\end{aligned}
\]

Integrability conditions \emph{($\der^2\equiv 0$)} 
of these equations imply the existence of a 1-form $\varpi_5$ such that:
\[
\begin{aligned}
\der\varpi_1
=&\,
\bar{\om}^1\dz\varpi_5+\bar{\om}^2\dz\varpi_4-\bar{\om}^4\dz\varpi_3,\\
\der \varpi_2
=&\,
-\tfrac14 {\color{blue}I^3}\bar{\om}^1\dz\varpi_3-\tfrac18I^3{}_{|5}\bar{\om}^1\dz\varpi_4-\tfrac12\bar{\om}^2\dz\varpi_4-\tfrac12\bar{\om}^4\dz\varpi_3+\\
&\ \ \ 
\tfrac{1}{16}(I^3{}_{|522}+2I^3{}_{|42}-8I^2{}_{|5})\bar{\om}^1\dz\bar{\om}^2+\tfrac{1}{16}(I^3{}_{|523}+2I^3{}_{|43})\bar{\om}^1\dz\bar{\om}^3+\\
&\ \ \ 
\tfrac{1}{16}(8I^1{}_{|3}-I^3{}_{|524}-2I^3{}_{|44})\bar{\om}^1\dz\bar{\om}^4-\tfrac{1}{16}(I^3{}_{|525}+2I^3{}_{|45})\bar{\om}^1\dz\bar{\om}^5+\\
&\ \ \
\tfrac{1}{8}(I^3{}_{|52}-2I^3{}_{|4})\bar{\om}^2\dz\bar{\om}^4-\tfrac12I^3{}_{|5}\bar{\om}^2\dz\bar{\om}^5+{\color{blue}I^3}\bar{\om}^3\dz\bar{\om}^4-\bar{\om}^3\dz\bar{\om}^5,\\
\der \varpi_3
=&\,
\varpi_3\dz(\tfrac12\varpi_1+\varpi_2)+\tfrac18(2I^3{}_{|4}+I^3{}_{|52})\bar{\om}^1\dz\varpi_3+\tfrac14 {\color{blue}I^3}\bar{\om}^2\dz\varpi_3+\\
&\ \ \
\tfrac18I^3{}_{|5}\bar{\om}^2\dz\varpi_4+\tfrac12\bar{\om}^2\dz\varpi_5+\bar{\om}^3\dz\varpi_4+J^1\bar{\om}^1\dz\bar{\om}^2+\\
&\ \ \
\tfrac{1}{4}(4I^2{}_{|5}+4I^3{}_{|1}-2I^3{}_{|42}-I^3{}_{|522})\bar{\om}^1\dz\bar{\om}^3+({\color{red}I^1}{\color{blue}I^3}-I^1{}_{|2})\bar{\om}^1\dz\bar{\om}^4+\\
&\ \ \ 
{\color{red}I^1}\bar{\om}^1\dz\bar{\om}^5-\tfrac{1}{16}(2I^3{}_{|43}+I^3{}_{|523})\bar{\om}^2\dz\bar{\om}^3+\tfrac{1}{16}(I^3{}_{|524}-8I^1{}_{|3}+2I^3{}_{|44})\bar{\om}^2\dz\bar{\om}^4+
\\
&\ \ \
\tfrac{1}{16}(2I^3{}_{|45}+I^3{}_{|525})\bar{\om}^2\dz\bar{\om}^5-\tfrac18(2I^3{}_{|4}+I^3{}_{|52})\bar{\om}^3\dz\bar{\om}^4,\\
\der \varpi_4
=&\,
\varpi_4\dz(\tfrac12\varpi_1-\varpi_2)+\tfrac18(2I^3{}_{|4}+I^3{}_{|52})\bar{\om}^1\dz\varpi_4-\tfrac14 {\color{blue}I^3}\bar{\om}^4\dz\varpi_3-\\
&\ \ \
\tfrac18I^3{}_{|5}\bar{\om}^4\dz\varpi_4+\tfrac12\bar{\om}^4\dz\varpi_5+\bar{\om}^5\dz\varpi_3+
\\
&\ \ \
\tfrac{1}{128}\Big(16(I^3{}_{|14}-{\color{red}I^1}I^3{}_{|3})+8(I^3{}_{|521}-I^1{}_{|3}{\color{blue}I^3})+2{\color{blue}I^3}I^3{}_{|44}+{\color{blue}I^3}I^3{}_{|524}\Big)\bar{\om}^1\dz\bar{\om}^4+\\
&\ \ \
\tfrac{1}{2}(2I^2{}_{|4}+{\color{green}I^2}I^3{}_{|5})\bar{\om}^1\dz\bar{\om}^2-{\color{green}I^2}\bar{\om}^1\dz\bar{\om}^3+\tfrac{1}{16}(8I^2{}_{|5}-2I^3{}_{|42}-I^3{}_{|522})\bar{\om}^2\dz\bar{\om}^4+\\
&\ \ \
\tfrac{1}{4}(I^3{}_{|524}-4I^1{}_{|3}+2I^3{}_{|44}+2I^3{}_{|51})\bar{\om}^1\dz\bar{\om}^5+\tfrac{1}{8}(2I^3{}_{|4}+I^3{}_{|52})\bar{\om}^2\dz\bar{\om}^5-
\\
&\ \ \
\tfrac{1}{16}(2I^3{}_{|43}+I^3{}_{|523})\bar{\om}^3\dz\bar{\om}^4-\tfrac{1}{16}(2I^3{}_{|45}+I^3{}_{|525})\bar{\om}^4\dz\bar{\om}^5\end{aligned}
\]
\[
\begin{aligned}
\der \varpi_5
=&\,
\varpi_5\dz\varpi_1+2\varpi_4\dz\varpi_3+J^2\bar{\om}^1\dz\varpi_3+J^3\bar{\om}^1\dz\varpi_4+\tfrac14(2I^3{}_{|4}+I^3{}_{|52})\bar{\om}^1\dz\varpi_5+\\
&\ \ \
\tfrac18(2I^3{}_{|4}+I^3{}_{|52})\bar{\om}^2\dz\varpi_4-\tfrac18(2I^3{}_{|4}+I^3{}_{|52})\bar{\om}^4\dz\varpi_4+J^4\bar{\om}^1\dz\bar{\om}^2+J^5\bar{\om}^1\dz\bar{\om}^3+
\\
&\ \ \
J^6 \bar{\om}^1\dz\bar{\om}^4+J^7\bar{\om}^1\dz\bar{\om}^5-{\color{green}I^2}\bar{\om}^2\dz\bar{\om}^3+J^8\bar{\om}^2\dz\bar{\om}^4+
\\
&\ \ \
\tfrac{1}{4}(I^3{}_{|524}-4I^1{}_{|3}+2I^3{}_{|44}+2I^3{}_{|51})\bar{\om}^2\dz\bar{\om}^5+
\\
&\ \ \
\tfrac{1}{4}(4I^2{}_{|5}+4I^3{}_{|1}-2I^3{}_{|42}-I^3{}_{|522})\bar{\om}^3\dz\bar{\om}^4-{\color{red}I^1}\bar{\om}^4\dz\bar{\om}^5.
\end{aligned}
\]
\[
\begin{aligned}\der I^3{}_{|2}
=&\,
\tfrac{1}{16}\Big(16(I^3{}_{|12}-{\color{green}I^2}I^3{}_{|5})+{\color{blue}I^3}(8I^2{}_{|5}-2I^3{}_{|42}-I^3{}_{|522})\Big)\bar{\om}^1+I^3{}_{|22}\bar{\om}^2+I^3{}_{|23}\bar{\om}^3+
\\
&\ \ \ 
\tfrac18\Big(8(I^3{}_{|42}+I^3{}_{|1})+{\color{blue}I^3}(I^3{}_{|52}-2I^3{}_{|4})\Big)\bar{\om}^4+\tfrac12\Big(2(I^3{}_{|52}-I^3{}_{|4})-{\color{blue}I^3}I^3{}_{|5}\Big)\bar{\om}^5-\\&
I^3{}_{|2}\varpi_1+2I^3{}_{|2}\varpi_1-I^3{}_{|3}\varpi_3-{\color{blue}I^3}\varpi_4,
\\
\der I^3{}_{|3}
=&\,
\tfrac{1}{16}\Big(16I^3{}_{|13}-2I^3{}_{|3}(2I^3{}_{|4}+I^3{}_{|52})-{\color{blue}I^3}(I^3{}_{|523}+2I^3{}_{|43})\Big)\bar{\om}^1+(I^3{}_{|23}-{\color{blue}I^3}I^3{}_{|3})\bar{\om}^2+
\\
&\ \ \
I^3{}_{|33}\bar{\om}^3+
\tfrac12\Big(I^3{}_{|523}+2I^3{}_{|43}-2(I^3{}_{|2}+({\color{blue}I^3})^2)\Big)\bar{\om}^4+3{\color{blue}I^3}\bar{\om}^5-
\tfrac12I^3{}_{|3}\varpi_1+3I^3{}_{|3}\varpi_2,
\\
\der I^3{}_{|5}
=&\,
I^3{}_{|51}\bar{\om}^1+I^3{}_{|52}\bar{\om}^2+4{\color{blue}I^3}\bar{\om}^3+I^3{}_{|54}\bar{\om}^4+I^3{}_{|55}\bar{\om}^5-\tfrac12I^3{}_{|5}\varpi_1-I^3{}_{|5}\varpi_2,\\
\der I^3{}_{|52}=
&\,
I^3{}_{|521}\bar{\om}^1+I^3{}_{|522}\bar{\om}^2+4\big(({\color{blue}I^3})^2+I^3{}_{|2}\big)\bar{\om}^3+I^3{}_{|524}\bar{\om}^4+
\\
&\ \ \
\big(2I^3{}_{|45}+(I^3{}_{|5})^2+I^3{}_{|525}-2I^3{}_{|54}\big)\bar{\om}^5-I^3{}_{|52}\varpi_1-4{\color{blue}I^3}\varpi_3,
\\
\der I^3{}_{|55}
=&\,
\tfrac{1}{16}\big(16I^3{}_{|515}-I^3{}_{|5}(2I^3{}_{|45}+I^3{}_{|525})-2I^3{}_{|55}(I^3{}_{|4}+I^3{}_{|52})\big)\bar{\om}^1+\\
&\ \ \
\tfrac12\big(4I^3{}_{|45}+(I^3{}_{|5})^2+2I^3{}_{|525}-2I^3{}_{|54}\big)\bar{\om}^2+3I^3{}_{|5}\bar{\om}^3+\tfrac12(2I^3{}_{|545}+I^3{}_{|5}I^3{}_{|55})\bar{\om}^4+
\\
&\ \ \
I^3{}_{|555}\bar{\om}^5-\tfrac12I^3{}_{|55}\varpi_1-3I^3{}_{|55}\varpi_2.
\end{aligned}
\]
Here, the coefficients ${\color{red}I^1}$, ${\color{green}I^2}$,
${\color{blue}I^3}$ are the respective incarnations of the basic
para-CR relative invariants ${\color{red}\inc{A}{}}$,
${\color{green}\inc{B}{}}$ and ${\color{blue}\inc{C}{}}$ from
Theorem~{\ref{the1}}. Each of them is a nonzero multiple of the
respective ${\color{red}\inc{A}{}}$, ${\color{green}\inc{B}{}}$,
${\color{blue}\inc{C}{}}$, as follows:
\[
\begin{aligned}
{\color{red}I^1}\sim&~{\color{red}9D^2H_r-27DH_p-18H_rDH_r+18H_pH_r+4H_r^3+54H_z},\\
{\color{green}I^2}\sim&~{\color{green}40G_{ppp}^3-45G_{pp}G_{ppp}G_{pppp}+9G_{pp}^2G_{ppppp}},\\
{\color{blue}I^3}\sim&~{\color{blue}2G_{ppp}+G_{pp}H_{rr}}.
\end{aligned}
\]
The other functions, such as e.g. $I^3{}_{|5}$, are coframe derivatives of the basic invariants ${\color{red}I^1}$,  ${\color{green}I^2}$ and ${\color{blue}I^3}$, with the convention that, for a function 
$f$:
\[
\der f
=
f_{|1}\om^1+f_{|2}\om^2+f_{|3}\om_3+f_{|4}\om^4+f_{|5}\om^5+(\dots)\varpi_1+(\dots)\varpi_2+(\dots)\varpi_3+(\dots)\varpi_4.
\]
The dotted coeffcients in this expression follow from $\der^2=0$
applied to the above EDS and to $f$. The coefficients
$J^1,J^2,\dots,J^8$ are not important here.
%\begin{aligned}128J^1=&64(2I^1{}_{|23}+2I^2{}_{|45}-I^3{}_{|424}+I^2{}_{|5}I^3{}_{|5})-120I^1{}_{|3}I^3+112(I^3{}_{|14}-I^3{}_{|3}I^1)-\\&2I^3{}_{|44}I^3-8I^3{}_{|521}-32I^3{}_{|5224}-I^3{}_{|524}I^3 \end{aligned}$$
\end{theorem}

In this section, we have an {\em a priori} `crazy idea' of relating
the EDS of Theorem~{\ref{susend}} to the EDS~{\eqref{e.c.dtheta_10d}}
from Theorem~{\ref{th.c.1}} describing 3\textsuperscript{rd} order
ODEs. There are several reasons indicating that this idea is not so
weird as it looks at first glance.

\begin{itemize}

\smallskip\item 
As we already noticed, in our 5-dimensional para-CR structure theory,
there is a family of third order ODEs $z_{xxx}=H(x,y,z,z_x,z_{xx})$
incorporated.

\smallskip\item 
One of our para-CR invariants ${\color{red}\inc{A}{}}$ is the contact
(therefore also point) W\"unschmann invariant appearing in the theory
of 3\textsuperscript{rd} order ODEs.

\smallskip\item 
The flat model of our 5-dimensional para-CR structures is described in
terms of the Maurer-Cartan forms on the Lie group $Sp(4,\bbR)$, which
is the same as the description of the flat model for the geometry of
third order ODEs considered modulo contact transformation of
variables, which is also given as an EDS satisfied by the
Maurer-Cartan forms on $Sp(4,\bbR)$.
\end{itemize}

%The more convincing argument for our `crazy idea' is as follows:

%Given a 3\textsuperscript{rd} order ODE, $z'''=H(x,z,z',z'')$, one encodes its geometric properties in the coframe
%$$\begin{aligned}
%  \alpha^1=&\der z-p\der x,\\
%  \alpha^2=&\der p-r\der x,\\
%  \alpha^3=&\der r-H\der x,\\
%  \alpha^4=&\der x,\end{aligned}$$
%where we abreviated $z'=p$, $z''=r$.
%Then, if the ODE is considered modulo contact transformations of variables, these forms are defined up to the following transformations $(\alpha^1,\alpha^2,\alpha^3,\alpha^4)\to (\bar{\alpha}^1,\bar{\alpha}^2,\bar{\alpha}^3,\bar{\alpha}^4)$ with:
%$$\begin{aligned}
%  \bar{\alpha}^1=f_1\alpha^1\\
%  \bar{\alpha}^2=f_2\alpha^1+\rho\mathrm{e}^\phi\alpha^2\\
%  \bar{\alpha}^3=f_5\alpha^1+f_6\alpha^2+f_7\alpha^3\\
%  \bar{\alpha}^4=\bar{f}_2\alpha^1+f_8\alpha^3+f_9\alpha^4.
% \end{aligned}$$

This motivates our `crazy question', which actually, due to the
discrete symmetry $D_1 \longleftrightarrow D_2$ between the two
integrable para-CR distributions $D_1$ and $D_2$, consists of two
questions:

\medskip\noindent{\bf Q1.}
{\em Can we bring the EDS of Theorem~{\ref{susend}}, by only using
para-CR transformations of forms $(\bar{\omega}^1, \bar{\omega}^2,
\bar{\omega}^3, \bar{\omega}^4, \bar{\omega}^5)$, to the
EDS~{\eqref{e.c.dtheta_10d}} describing contact equivalence classes
of 3\textsuperscript{rd} order ODEs?  More specifically, can we
force the system of 1-forms}
\begin{equation}
\begin{aligned}
\theta^1
=&\,
f_1\bar{\om}^1,\\
\theta^2
=&\,
f_2\bar{\om}^1+\rho\mathrm{e}^\phi\bar{\om}^2+f_4\bar{\om}^3,\\
\theta^3
=&\,
f_5\bar{\om}^1+f_6 \bar{\om}^2+f_7\bar{\om}^3,\\
\theta^4
=&\,
\bar{f}_2\bar{\om}^1+\rho\mathrm{e}^{-\phi}\bar{\om}^4+\bar{f}_4\bar{\om}^5,\\
\theta^5
=&\,
\bar{f}_5\bar{\om}^1+\bar{f}_6\bar{\om}^4+\bar{f}_7\bar{\om}^5  ,\end{aligned}\label{godwun}
\end{equation}
{\em to satisfy the EDS~{\eqref{e.c.dtheta_10d}},  
by an appropriate choice of the fiber parameters $(f_1, f_2, \rho, \phi, f_4, f_5, f_6, f_7, \bar{f}_2, \bar{f}_4, \bar{f}_5, \bar{f}_6,
\bar{f}_7)$?}

\medskip\noindent{\bf Q2.}
{\em The same question as {\bf Q1}, but now with the flip
$(\bar{\om}^2, \bar{\om}^3) \longleftrightarrow
(\bar{\om}^4, \bar{\om}^5)$,
namely: can we force the system of 1-forms}
\begin{equation}
\begin{aligned}
\theta^1
=&\,
f_1\bar{\om}^1,\\
\theta^2
=&\,
f_2\bar{\om}^1+\rho\mathrm{e}^\phi\bar{\om}^4+f_4\bar{\om}^5,\\
\theta^3
=&\,
f_5\bar{\om}^1+f_6 \bar{\om}^4+f_7\bar{\om}^5,\\
\theta^4
=&\,
\bar{f}_2\bar{\om}^1+\rho\mathrm{e}^{-\phi}\bar{\om}^2+\bar{f}_4\bar{\om}^3,\\
\theta^5
=&\,
\bar{f}_5\bar{\om}^1+\bar{f}_6\bar{\om}^2+\bar{f}_7\bar{\om}^3,
\end{aligned}\label{godwn}
\end{equation}
{\em to satisfy the EDS~{\eqref{e.c.dtheta_10d}},  
by an appropriate choice of the fiber parameters $(f_1, f_2, \rho, \phi, f_4, f_5, f_6, f_7, \bar{f}_2, \bar{f}_4, \bar{f}_5, \bar{f}_6, 
\bar{f}_7)$?}

\medskip

The next theorem gives the {\em if-and-only-if} answer for these
questions, as well as the obstructions to achive the goals specified
in questions {\bf Q1} and {\bf Q2}, in terms of the para-CR
invariants.

\begin{theorem}\label{bril}$\:$

$\bullet$\,
Question {\bf Q1} above has a positive answer if and only if
$I^3{}_{|3}\equiv 0$. The para-CR structures related to
$I^3{}_{|3}\equiv 0$ can be distinguished by the $\spa(4,\bbR)$-valued
Cartan connection~{\eqref{e.c.conn_sp}} whose curvature $\hat{K}$ has
the basic invariant ${\color{blue}\inc{Z}{1}\equiv 0}$ and the basic
invariant ${\color{red}\inc{A}{1}}$ proportional to the W\"unschmann
invariant
\[
{\color{red}\inc{A}{1}}~\sim~{\color{red}9D^2H_r-27DH_p-18H_rDH_r+18H_pH_r+4H_r^3+54H_z}.
\]

$\bullet$\,
Question {\bf Q2} above has a positive answer if and only if
$I^3{}_{|55} \equiv 0$. The para-CR structures related to $I^3{}_{|55}
\equiv 0$ can be distinguished by the $\spa(4, \bbR)$-valued Cartan
connection~{\eqref{e.c.conn_sp}} whose curvature $\hat{K}$ has the
basic invariant ${\color{blue} \inc{Z}{1} \equiv 0}$ and the basic
invariant ${\color{red} \inc{A}{1}}$ proportional to the W\"unschmann
invariant
\[
{\color{green}\inc{A}{1}}~\sim~{\color{green}40 G_{ppp}^3-45G_{pp}G_{ppp}G_{pppp}+9G_{pp}^2G_{ppppp}}.
\]

$\bullet$\,
Furthermore, each condition $I^3{}_{|3}\equiv 0$ and
$I^3{}_{|55}\equiv 0$, considered separately, implies that the
relative fundamental para-CR invariant ${\color{blue}I^3}\equiv 0$. So
there is only \emph{one} `if and only if' condition for a positive
answer to questions {\bf Q1} or {\bf Q2}: any of them has a positive
answer if and only if the para-CR invariant ${\color{blue}\inc{C}{}}$
vanishes:
\[
{\color{blue}2G_{ppp}+G_{pp}H_{rr}\equiv 0}.
\]
\end{theorem}

\begin{remark}\label{rek}
Before starting the proof, we remark that this theorem provides a way
of \emph{transforming two classical invariants, the W\"unschmann one
and the Monge one, into each other}. This can be achieved by passing
from the third order ODE corresponding to the Cartan connection
related to the question {\bf Q1}, to its \emph{dual}
3\textsuperscript{rd} order ODE, described by the Cartan connection
related to question {\bf Q2}.    
\end{remark}

\proof[Proof of Theorem~{\ref{bril}}]
%/home/pawel/niebieski/notebooks/merker/homogeneous models/da_capo_2019/december_2019/ala_godlinski_wunschman.nb
We first answer question {\bf Q1}. We start with the forms
$(\theta^1,\theta^2,\theta^3,\theta^4,\theta^5)$ as
in~{\eqref{godwun}}, and we try to make normalizations on
$\der\theta^i$ as in~{\eqref{e.c.dtheta_10d}}. For full generality, we
will not use~{\eqref{e.c.dtheta_10d}} with the 1-form $\theta^5$ in
it. We will call this 1-form $\Omega_0$ for a while. As we will see in
the proof, the procedure we apply now, which is an adaptation of
Cartan's equivalence method, is powerfull enough to determine the
relation between $\Omega_0$ and $theta^5$.

The first normalizations coming from~{\eqref{e.c.dtheta_10d}}, namely
$\der\theta^1 \w \theta^1 \w \theta^2=0$ and
$\der\theta^1 \w \theta^1 \w \theta^4=0$, give
\[
f_4=\bar{f}_4=0,
\]
and then, $\der\theta^1\w\theta^1 = -\theta^1\w\theta^2\w\theta^4$, 
gives 
\[
f_1=-\rho^2.
\]
Now the first condition in~{\eqref{e.c.dtheta_10d}} 
enables to determine 
\[
\Omega_1
=
\varpi_1-\frac{\bar{f}_2}{\rho^2}\theta^2+\frac{f_2}{\rho^2}\theta^4+\der\log(\rho^2)-u_1\theta^1,
\]
up to the term with $\theta^1$, which requires to introduce a new
variable $u_1$.

We now make the normalization $\der\theta^2\w\theta^1\w\theta^2=-\theta^1\w\theta^2\w\theta^3\w\theta^4,$ which results in 
\[
f_7
=
-\mathrm{e}^{2\phi}.
\]
After this normalization, the second equation
in~{\eqref{e.c.dtheta_10d}} solves for $\Omega_2$ and $\Omega_3$ as
follows:
\[
\begin{aligned}
\Omega_2
=&\,
-\frac{\bar{f}_2}{\rho^2}\theta^3+\frac{f_5\rho^2-f_2^2-f_2f_6\rho\mathrm{e}^{-\phi}}{\rho^4}\theta^4-\frac{f_2}{\rho^2}(\tfrac12\varpi_1+\varpi_2)-\frac{\mathrm{e}^\phi}{\rho}\varpi_3+\frac{f_2}{\rho^2}\der\log(\frac{\rho\mathrm{e}^\phi}{f_2})+\\
&\ \ \
\frac{2\rho^2u_2-f_2u_1}{2\rho^2}\theta^1+\frac{2\rho^4u_3-f_2\bar{f}_2}{2\rho^4}\theta^2,\\
\Omega_3
=&\,
-\frac{f_2+f_6\rho\mathrm{e}^{-\phi}}{\rho^2}\theta^4+\tfrac12\varpi_1-\varpi_2+\der\log(\rho\mathrm{e}^\phi)+\frac{2\rho^4u_3-3f_2\bar{f}_2-2\bar{f}_2f_6\rho\mathrm{e}^{-\phi}}{2\rho^4}\theta^1-\frac{\bar{f}_2+2\rho^2u_4}{2\rho^2}\theta^2,
\end{aligned}
\]
where $u_2,u_3,u_4$ are new variables taking account on how
indeterminate are $\Omega_2$ and $\Omega_3$.

Now $\der\theta^4 \w \theta^2 \w \theta^4 = \Omega_4 \w \theta^1 \w
\theta^2\w\theta^4$ gives
\[
\begin{aligned}
\Omega_4=&\frac{f_2\mathrm{e}^{-2\phi}}{\bar{f}_7\rho^2}\theta^5+\frac{\bar{f}_2}{\rho^2}(\varpi_2-\tfrac12\varpi_1)-\frac{\mathrm{e}^{-\phi}}{\rho}\varpi_4+\frac{\bar{f}_2}{\rho^2}\der\log(\frac{\rho\mathrm{e}^{-\phi}}{\bar{f}_2})+\frac{3f_2\bar{f}_2^2-2\bar{f}_2\rho^4(u_1+u_3)-2\rho^6u_5+2\bar{f}_2^2f_6\rho\mathrm{e}^{-\phi}}{2\rho^6}\theta^1+\\
&\frac{2\bar{f}_2\rho^2u_4-\bar{f}_2^2-2\rho^4u_6}{2\rho^4}\theta^2+\frac{2f_2\bar{f}_2-\rho^4u_7+\bar{f}_2f_6\rho\mathrm{e}^{-\phi}}{\rho^4}\theta^4,
\end{aligned}
\]
and $\der\theta^4 \w \theta^4 = (\Omega_4\w\theta^1 +
\Omega_0\w\theta^2)\w\theta^4,$ shows that the ODE 1-form $\Omega_0$
must be expressed in terms of $\theta^1$, $\theta^2$, $\theta^4$,
$\theta^5$ as follows:
\[
\begin{aligned}
\Omega_0=\frac{\mathrm{e}^{-2\phi}}{\bar{f}_7}\theta^5+\frac{2\bar{f}_7\rho^2(\bar{f}_2 u_4-\rho^2u_6)-3\bar{f}_2^2\bar{f}_7-2\bar{f}_2\bar{f}_6\rho\mathrm{e}^{-\phi}+2\bar{f}_5\rho^2\mathrm{e}^{-2\phi}}{2\bar{f}_7\rho^4}\theta^1+u_8\theta^2+u_9\theta^4.
\end{aligned}
\]
All of this is true with new undetermined variables $u_5,\ldots,u_9$.
Please note that in this formula there are \emph{no terms consistsing
of the differentials of the group parameters}!

Now to get the fourth equation~{\eqref{e.c.dtheta_10d}} satisfied we
have to put
\[
u_7=\frac{3f_2\bar{f}_2\bar{f}_7+2\bar{f}_7\rho^4(u_1+u_3)+2f_2\bar{f}_6\rho\mathrm{e}^{-\phi}}{2\bar{f}_7\rho^4}\quad\mathrm{and}\quad u_9=\frac{2\bar{f}_7\rho^2u_4-3\bar{f}_2\bar{f}_7-2\bar{f}_6\rho\mathrm{e}^{-\phi}}{2\bar{f}_7\rho^2}.
\]
With these normalizations and definitions of $\Omega_\mu$'s, the
differentials $\der\theta^1$, $\der\theta^2$, $\der\theta^4$ are
precisely as in~{\eqref{e.c.dtheta_10d}}.

The third equation in~{\eqref{e.c.dtheta_10d}} is achieved by a unique
choice of $f_5$, $f_6$, $u_4$ and $u_3$ as follows:
\[
f_5=-\frac{f_2^2}{2\rho^2},\quad 
f_6=-\frac{f_2\mathrm{e}^\phi}{\rho},\quad 
u_4=\frac{\bar{f}_2-{\color{blue}I^3}\rho\mathrm{e}^{-\phi}}{2\rho^2},\quad 
u_3=\frac{8f_2 {\color{blue}I^3}\mathrm{e}^{-\phi}-\rho(2I^3{}_{|4}+I^3{}_{|52}+8u_1\rho^2)}{16\rho^3}.
\]

With these normalizations, we achieve that $\der\theta^3$ is precisely
as in the third equation in~{\eqref{e.c.dtheta_10d}} with
\[
{\color{red}\inc{A}{1}}
=
-(\frac{\mathrm{e}^\phi}{\rho})^3{\color{red}~I^1}
\quad\mathrm{and}\quad
\inc{A}{2}
=
\frac{f_2\big(4f_2\bar{f}_2+\rho^2(2I^3{}_{|4}+I^3{}_{|52})-4\rho^4u_1\big)+8\rho^6u_2-4f_2^2\rho\mathrm{e}^{-\phi} {\color{blue}I^3}}{8\rho^6}.
\]

Now there is a unique way to bring the differential $\der\Omega_0$ to
the form of the ninth equation~{\eqref{e.c.dtheta_10d}}. For this we
have:
\[
u_8=\frac{I^3{}_{|3}\mathrm{e}^{-3\phi}}{2\rho}, 
\quad\mathrm{and}\quad 
u_6=\frac{-4\bar{f}_2^2+\rho^2\mathrm{e}^{-2\phi}(2I^3{}_{|43}+I^3{}_{|523})-4f_2\rho\mathrm{e}^{-3\phi}I^3{}_{|3}}{8\rho^4}.
\]
After this normalization the formula for $\der\Omega_0$ is as
in~{\eqref{e.c.dtheta_10d}}. In particular, we get explicit
expressions for $\inc{A}{5}$, $\inc{Z}{4}$, $\inc{A}{7}+\inc{Z}{3}$,
which are not important here, but also the formula for
${\color{blue}\inc{Z}{1}}$ which is:
\[
{\color{blue}\inc{Z}{1}}
=
\frac{\mathrm{e}^{-5\phi}}{2\rho}~I^3{}_{|33}.
\]

The last 3\textsuperscript{rd} order ODE invariant 1-form $\Omega_5$
is now determined from $\der\Omega_1\w\theta^2\w\theta^4 =
\Omega_5\w\theta^1\w\theta^2\w\theta^4$, as
\[
\Omega_5
=
-u_1\Omega_1-\frac{\bar{f}_2}{\rho^2}\Omega_2+\frac{f_2}{\rho^2}\Omega_4+\frac{\bar{f}_2\mathrm{e}^\phi}{\rho^3}\varpi_3-\frac{f_2\mathrm{e}^{-\phi}}{\rho^3}\varpi_4+\frac{1}{\rho^2}\varpi_5-u_{10}\theta^1-u_{11}\theta^2-u_{12}\theta^4-\der u_1
\]
up to $\theta^1$, $\theta^2$, $\theta^4$ terms, which require
introduction of new parameters $u_{10}$, $u_{11}$, $u_{12}$.

Now there is a unique way of killing all the unwanted terms in
$\der\Omega_\mu$, $\mu=0,\ldots, 5$, to achieve the full
system~{\eqref{e.c.dtheta_10d}}. It turns out that now, this involves
solving \emph{linear} equations for all the remaining auxiliary
variables $u_2$, $u_5$, $u_{10}$, $u_{11}$, $u_{12}$\,\,---\,\,except
$u_1$. They are determined successively, as follows: $u_5$ is
determined by killing the unwanted terms in $\der\Omega_1\w\theta^4$,
$u_2$ is determined by killing the unwanted terms in $\der\Omega_1$,
$u_{12}$ is deteremined by killing the unwanted terms in
$\der\Omega_2\w\theta^1$, $u_{11}$ is determined by killing the
unwanted terms in $\der\Omega_2$, and finally $u_{10}$ is deteremined
by killing the unwanted terms in
$\der\Omega_5\w\theta^1\w\theta^3$. The explicit expressions for these
auxiliary variables are not relevant here.

The final result of these normalizations is:
 \[
\begin{aligned}
\theta^1=&-\rho^2\bar{\om}^1,\\
\theta^2=&f_2\bar{\om}^1+\rho\mathrm{e}^\phi\bar{\om}^2,\\
\theta^3=&-\frac{f_2^2}{2\rho^2}\bar{\om}^1-\frac{f_2\mathrm{e}^\phi}{\rho}\bar{\om}^2-\mathrm{e}^{2\phi}\bar{\om}^3,\\
\theta^4=&\bar{f}_2\bar{\om}^1+\rho\mathrm{e}^{-\phi}\bar{\om}^4,\\
\theta^5=&\bar{f}_5\bar{\om}^1+\bar{f}_6\bar{\om}^4+\bar{f}_7\om^5,
\end{aligned}
\]
\[
\begin{aligned}
\boxed{\Omega_0}=& \,\boxed{s_4\theta^1-\frac{2\bar{f}_2\bar{f}_7+\rho\mathrm{e}^{-\phi}(2\bar{f}_6+\bar{f}_7{\color{blue}I^3})}{2\bar{f}_7\rho^2}\theta^4+\frac{\mathrm{e}^{-2\phi}}{\bar{f}_7}\theta^5}+\frac{\mathrm{e}^{-3\phi}}{2\rho}~I^3{}_{|3}~\theta^2,\\~&\\
\Omega_1=&-u_1\theta^1-\frac{\bar{f}_2}{\rho^2}\theta^2+\frac{f_2}{\rho^2}\theta^4+\varpi_1+\der\log(\rho^2), \\~&\\
\Omega_2=&\frac{8I^1{}_{|3}\rho^3\mathrm{e}^\phi-3f_2\big(4f_2\bar{f}_2+\rho^2(2I^3{}_{|4}+I^3{}_{|52})\big)+12f_2^2\rho\mathrm{e}^{-\phi}{\color{blue}I^3}}{24\rho^6}\theta^1+\\&\frac{8f_2\rho\mathrm{e}^{-\phi}{\color{blue}I^3}-8f_2\bar{f}_2-\rho^2(2I^3{}_{|4}+I^3{}_{|52})-8\rho^4u_1}{16\rho^4}\theta^2-\frac{\bar{f}_2}{\rho^2}\theta^3-\frac{f_2^2}{2\rho^4}\theta^4-\frac{f_2}{\rho^2}(\tfrac12\varpi_1+\varpi_2)-\frac{\mathrm{e}^\phi}{\rho}\varpi_3+\\&\frac{f_2}{\rho^2}\der\log(\frac{\rho\mathrm{e}^\phi}{f_2}),\\~&\\
\Omega_3=&\frac{8f_2\rho\mathrm{e}^{-\phi}{\color{blue}I^3}-8f_2\bar{f}_2-\rho^2(2I^3{}_{|4}+I^3{}_{|52})-8\rho^4u_1}{16\rho^4}\theta^1+\frac{\rho\mathrm{e}^{-\phi}{\color{blue}I^3}-2\bar{f}_2}{2\rho^2}\theta^2-\frac{f_2+f_6\rho\mathrm{e}^{-\phi}}{\rho^2}\theta^4+\tfrac12\varpi_1-\varpi_2+\\&\der\log(\rho\mathrm{e}^\phi), \\~\\
\Omega_4=&s_1\theta^1+s_2\theta^2+s_3\theta^4+\frac{f_2\mathrm{e}^{-2\phi}}{\bar{f}_7\rho^2}\theta^5+\frac{\bar{f}_2}{\rho^2}(\varpi_2-\tfrac12\varpi_1)-\frac{\mathrm{e}^{-\phi}}{\rho}\varpi_4+\frac{\bar{f}_2}{\rho^2}\der\log(\frac{\rho\mathrm{e}^{-\phi}}{\bar{f}_2}),
\end{aligned}
\]
\[
\begin{aligned}
\Omega_5=&s_5\theta^1+s_6\theta^2+\frac{\bar{f}_2^2}{\rho^4}\theta^3+s_7 \theta^4+\frac{f_2^2\mathrm{e}^{-2\phi}}{\bar{f}_7\rho^4}\theta^5-u_1\varpi_1+\frac{2f_2\bar{f}_2}{\rho^4}\varpi_2+\frac{2\bar{f}_2\mathrm{e}^\phi}{\rho^3}\varpi_3-\frac{2f_2\mathrm{e}^{-\phi}}{\rho^3}\varpi_4+\frac{1}{\rho^2}\varpi_5-\der u_1-\\&\frac{2u_1}{\rho}\der\rho+\frac{\bar{f}_2}{\rho^4}\der f_2-\frac{f_2}{\rho^4}\der\bar{f}_2-\frac{2f_2\bar{f}_2}{\rho^4}\der\phi.
\end{aligned}
\]
Here
\[
s_4=\frac{-4\bar{f}_2^2\bar{f}_7-4\bar{f}_2(2\bar{f}_6+\bar{f}_7{\color{blue}I^3})\rho\mathrm{e}^{-\phi}+\big(8\bar{f}_5-\bar{f}_7(2I^3{}_{|43}+I^3{}_{|523})\big)\rho^2\mathrm{e}^{-2\phi}+4f_2\bar{f}_7I^3{}_{|3}\rho\mathrm{e}^{-3\phi}}{8\bar{f}_7\rho^4},
\]
and the coefficients $s_1,s_2,s_3,s_5,s_6,s_7$ although explicitely
determined, are totally irrelevant for the sequel.

The above 1-forms $(\theta^1, \ldots, \theta^4, \Omega_0, \ldots,
\Omega_5)$ satisfy the 3\textsuperscript{rd} order ODE
system~{\eqref{e.c.dtheta_10d}} with
\[
{\color{red}\inc{A}{1}}
=
-(\frac{\mathrm{e}^\phi}{\rho})^3{\color{red}~I^1}\quad\mathrm{and}\quad{\color{blue}\inc{Z}{1}}=\frac{\mathrm{e}^{-5\phi}}{2\rho}~I^3{}_{|33}.
\]
Also all other coefficients $\inc{A}{i}$, $\inc{Z}{j}$ and
$\inc{C}{1}$ are totally determined, as for example $\inc{A}{2} =
\frac{\mathrm{e}^\phi}{ 3\rho^3}I^1{}_{|3}$, or $\inc{A}{5} =
-\frac{\mathrm{e}^{-\phi}}{ 6\rho^3}I^1{}_{|33}$, but they are not
that illuminating to quote them here.

This explicitly shows that \emph{every} 5-dimensional para-CR
structure, which has Levi form degenerate in one direction and which
is not locally a trivial extension of a 3-dimensional nondegenerate
para-CR structure, defines an invariant EDS for a contact equivalence
class of 3\textsuperscript{rd} ODEs for which the classical
W\"unschmann invariant is the para-CR invariant
${\color{red}\inc{A}{}}$.

A problem arises if this obtained EDS is para-CR invariant. At first
glance yes, but it is really \emph{not}.

The reason for that is that the form $\theta^5$ disappeared from the
description. It was replaced by the form $\Omega_0$. Looking at the
explicit form of $\Omega_0$, one observes that the forms
$(\theta^1,\theta^2,\theta^3,\theta^4,\theta^5)$ and
$(\theta^1,\theta^2,\theta^3,\theta^4,\Omega_0)$ are \emph{not}
para-CR equivalent, {\em because} 
the 1-form $\theta^2$ appears in the formula
relating $\Omega_0$ and $\theta^5$. 
But $\theta^2$ should \emph{not} be there!
Only the `boxed' part of this formula consists of some para-CR
transformation between $\theta^5$ and $\Omega_0$. The appearence of
the term $\frac{\mathrm{e}^{-3\phi}I^3{}_{|3}}{2\rho}\theta^2$ in this
formula breaks the para-CR equivalence.

There is only one way to restore the para-CR invariance of the
obtained EDS: one has to \emph{restrict to para-CR structures for
which}
\[
I^3{}_{|3}\equiv 0.
\]
In such a case one can use the remaining para-CR transformations to
achieve
\[
\Omega_0=\theta^5
\]
reducing all the auxiliary parameters from $(f_1, f_2, \rho, \phi,
f_4, f_5, f_6, f_7, \bar{f}_2, \bar{f}_4, \bar{f}_5, \bar{f}_6,
\bar{f}_7, u_1,\ldots, u_{12})$ to only five $(\rho,\phi,f_2,
\bar{f}_2,u_1)$, This makes the resulting EDS really 10-dimensional,
as it should be for it to describe a curvature of a Cartan
$\spa(4,\bbR)$-connection.

The proof of the answer to the question {\bf Q2} is essentially the
same as above. We start with the lifted coframe~{\eqref{godwn}}, and
impose the normalizations required by the
system~{\eqref{e.c.dtheta_10d}} in the same order as in the case of
question {\bf Q1}. We skip the details, reporting here the important
differences only. The first of them is that now the normalizations
result in $\der\theta^3$ as in~{\eqref{e.c.dtheta_10d}}, but with
\[
{\color{green}\inc{A}{1}}
=
-(\frac{\mathrm{e}^\phi}{\rho})^3{\color{green}I^2}.
\]
The next difference is that now, 
in the induced EDS~{\eqref{e.c.dtheta_10d}}, 
the coefficient ${\color{blue}\inc{Z}{1}}$ is
\[
{\color{blue}\inc{Z}{1}}
=
-\frac{\mathrm{e}^{-5\phi}}{4\rho}I^3{}_{|555}.
\]
As the last important difference we mention 
that now the 1-form $\Omega_0$
appearing in the induced EDS~{\eqref{e.c.dtheta_10d}} 
is related to $\theta^5$ via:
\[
\boxed{\Omega_0}
=
\boxed{s_4\theta^1+\frac{4\bar{f}_2\bar{f}_7-\rho\mathrm{e}^{-\phi}(4\bar{f}_6-\bar{f}_7I^3{}_{|5})}{4\bar{f}_7\rho^2}\theta^4+\frac{\mathrm{e}^{-2\phi}}{\bar{f}_7}\theta^5}-\frac{\mathrm{e}^{-3\phi}}{4\rho}~I^3{}_{|55}~\theta^2
\]
So now, the term $\frac{\mathrm{e}^{ -3\phi}I^3{}_{|55}}{ 4\rho}
\theta^2$ brakes the para-CR equivalence, and to answer the question
{\bf Q2} in positive, we are forced to restrict to para-CR
structures with
\[
I^3{}_{|55}\equiv 0.
\]
Consequently, if we assume that $I^3{}_{|55}\equiv 0$, we finally use the remaining para-CR transformations to achieve
\[
\Omega_0=\theta^5.
\]
This again reduces all the auxiliary parameters from $(f_1,f_2, \rho,
\phi, f_4,f_5,f_6,f_7, \bar{f}_2, \bar{f}_4, \bar{f}_5,\bar{f}_6,
\bar{f}_7, u_1, \ldots, u_{12})$ to only five
$(\rho,\phi,f_2, \bar{f}_2,u_1)$, and makes the resulting EDS the
curvature conditions of a Cartan $\spa(4,\bbR)$-connection in 10
dimensions.

As the final step of the proof, we remark that if we insert any of the
conditions $I^3{}_{|3}\equiv 0$ or $I^3{}_{|55}\equiv 0$ into the
EDS~{\eqref{sysend}} then its integrability conditions ($\der^2\equiv
0$) very quickly show that each of them is equivalent to
\[
{\color{blue}I^3\equiv 0}.
\]
For this, observe that if $I^3{}_{|3}\equiv 0$, then the equation for
$\der I^3{}_{|3}$ in Theorem~{\ref{susend}} gives immediately
${\color{blue} I^3} \equiv 0$. Likewise, if $I^3{}_{|55}\equiv 0$ then
the equation for $\der I^3{}_{|55}$ in Theorem~{\ref{susend}} gives
$I^3{}_{|5}\equiv 0$, and then the equation for $\der I^3{}_{|5}$ in
gives eventually ${\color{blue} I^3}\equiv 0$.

This finishes the proof of Theorem~{\ref{bril}}.
\endproof

Theorem~{\ref{susend}}, and the calculations done in its proof,
have an interesting 

\begin{corollary}\label{coco}
Consider a 5-dimensional para-CR structure given by the system of
PDEs~{\eqref{ss}} in terms of functions $H=H(x,y,z,p,r)$ and
$G=G(x,y,z,p,r)$ satisfying conditions~{\eqref{ic}}, {\eqref{1deg}},
{\eqref{2ndeg}}. Assume for this structure that the para-CR invariant
${\color{blue}\inc{C}{}}$ vanishes:
\[
{\color{blue}2G_{ppp}+G_{pp}H_{rr}
\equiv 
0}.
\]
Then, associated to such a para-CR structure, there are \emph{two}
contact equivalence classes of third order ODEs. Both of these classes
of ODEs have their respective Chern
invariants zero:
\[
{\color{blue}\inc{Z}{}\equiv 0}.
\]

The other basic contact invariant of these (contact invariant) classes
of ODEs, namely the W\"unschmann invariant $\inc{A}{1}$ is
proportional:

\smallskip\noindent{\bf a)}
to the W\"unschmann para-CR invariant,
$\inc{A}{1}~\sim~{\color{red}9D^2H_r - 27DH_p - 18H_rDH_r + 18H_pH_r +
4H_r^3 + 54H_z}$, in the case of the first class of ODEs; and:
  
\smallskip\noindent{\bf b)}
to the Monge para-CR invariant, $\inc{A}{1}~\sim~{\color{green}40
G_{ppp}^3 - 45G_{pp}G_{ppp}G_{pppp} + 9G_{pp}^2G_{ppppp}}$, in the
case of the second class of ODEs.
\end{corollary}

\proof[Proof of Corollary~{\ref{coco}}]
In the language of Theorem~{\ref{bril}}, the assumption that
${\color{blue}\inc{C}{}}\equiv 0$ means that ${\color{blue}I^3}\equiv
0$. This, in particular means that $I^3{}_{|33}\equiv 0$ and that
$I^3{}_{|555}\equiv 0$. Thus, the quantity ${\color{blue}\inc{Z}{1}}$
vanishes in the EDS obatained from the normalizations of the lifted
coframe~{\eqref{godwun}} as well as of the lifted
coframe~{\eqref{godwn}}.

Moreover, since $I^3\equiv 0$ implies also $I^3{}_{|3}\equiv 0$ and
$I^3{}_{|55}\equiv 0$, we know from Theorem~{\ref{bril}} that both
EDSs with ${\color{blue}\inc{Z}{1}}\equiv 0$ are para-CR
invariant. According to Chern's theory of 3\textsuperscript{rd} order
ODEs considered modulo contact transformation~{\cite{Chern-1940,
Godlinski-Nurowski-2009}}, both EDS's, considered separately,
describe a contact equivalence class of 3\textsuperscript{rd} order
ODEs on the 4-dimensional leaf space $\mathrm{J}^2$ of the rank 6
integrable distribution annihilating 1-forms $\theta^1, \theta^2,
\theta^3, \theta^4$. This space $\mathrm{J}^2$ can be locally
identified with the space of second jets of the corresponding class of
3\textsuperscript{rd} order ODEs. This class, in both EDSs, has
Chern invariant equal to zero (because ${\color{blue}\inc{C}{}}\equiv
0$ implies ${\color{blue}\inc{Z}{1}}\equiv 0$ in the EDSs), and as it
visible from the proof of Theorem~{\ref{bril}}, the classical
W\"unschmann invariant $\inc{A}{1}$ either proportional to
${\color{red}9D^2H_r-27DH_p-18H_rDH_r+18H_pH_r+4H_r^3+54H_z}$, or to
${\color{green}40 G_{ppp}^3 - 45G_{pp} G_{ppp}G_{pppp} +
9G_{pp}^2G_{ppppp}}$, depending which of the two EDSs we 
are considering.

This finishes the proof of Corollary~{\ref{coco}}. 
For further details about it,
consult our Appendix in Section~{\ref{api}}.
\endproof

\proof[End of proof of Theorem~{\ref{th0}}]
Since in both of the 10-dimensional para-CR invariant EDSs we have
${\color{blue}\inc{Z}{1}\equiv 0}$, then according to the result of
Godli\'nski~{\cite{Godlinski-2008,Godlinski-Nurowski-2009}}, the image
of the projection $\pi_1:\mathrm{J}^2\to \mathrm{J}^1$ from the second
jet space $\mathrm{J}^2$ appearing in the proof of the above
corollary, which can be identified with the 3-dimensional leaf space
of the rank 7 integrable distribution annihilating 1-forms $\theta^1,
\theta^2, \theta^4$, acquires a natural 3-dimensional contact
projective geometry. This proves our Theorem~{\ref{th0}} from the
Introduction.
\endproof

To illustrate the phenomena described in this section we consider the
following Example.
%ala_godlinski_wunschman_I3eq0_example_inny.nb  
%/example_transformation_3b_drugie.nb 

\begin{example} 
Our starting point in this example is a para-CR structure defined in
terms of PDEs
\begin{equation}
\boxed{z_y=f(z_x)\quad\&\quad z_{xxx}=-z_{xx}^2\frac{f^{(3)}(z_{x})}{f''(z_x)},\quad \mathrm{for}\quad z=z(x,y),}\label{cre}
\end{equation}
with $f=f(p)$ being a differentiable function such that $f''(p)\neq 0$. In other words we have
\begin{equation}
G=f(p)
\quad\&\quad 
H=-r^2\frac{f^{(3)}(p)}{f''(p)}.\label{crh}
\end{equation}
It is straightforward to check that $\triangle H=D^3G$, $G_{pp}\neq 0$
and, more importantly, that 
\[
{\color{blue}2G_{ppp}+G_{pp}H_{rr}
\equiv
0}.
\]

Therefore the Theorem{\ref{bril}} and Corolary~{\ref{coco}} apply, and
we should see two equivalence classes of 3\textsuperscript{rd} order
ODEs associated with these para-CR structures as well as two contact
projective structures on the spaces of first jets for these ODEs.

Before passing to show how these structure are explicitly visible
here, we calculate the W\"unschmann invariant ${\color{red}\inc{A}{}}$
for the function $H$ from \eqref{crh}. We have:
\[
\begin{aligned}
{\color{red}\inc{A}{}}=&~{\color{red}9D^2H_r-27DH_p-18H_rDH_r+18H_pH_r+4H_r^3+54H_z}=\frac{r^3}{f''{}^3}\,\Big(40f^{(3)}{}^3-45f''f^{(3)}f^{(4)}+9f''{}^2f^{(5)}\Big)\\
&=
\,(\frac{r^3}{G_{pp}^3})\,\Big({\color{green}40 G_{ppp}^3-45G_{pp}G_{ppp}G_{pppp}+9G_{pp}^2G_{ppppp}}\Big)=2r^3{\color{green}\inc{B}{}}.
\end{aligned}
\]
Thus, for our para-CR structure, represented by the functions $G$ and
$H$, the W\"unschmann invariant $\color{red}\inc{A}{}$ is a
nonvanishing multiple of the Monge invariant $\color{green}\inc{B}{}$.

This is a \emph{special case} of the phenomenon mentioned in
Remark~{\ref{rek}}: in this example we found the \emph{explicit}
transformation between the W\"unschmann invariant for $H$ and Monge
invariant of $G$. It was possible explicitly because this example is
so special that, as we see in a minute, the two \emph{a'priori
different} contact equivalent classes of 3\textsuperscript{rd} order
ODEs naturally associated to our para-CR structure, are actualy
\emph{the same}.

To see this we first write the coframe on $M^5$ encoding our para-CR
structure. This is given by:
\[
\begin{aligned}
\om^1=&\,\der z-p\der x-f\der y,\\
\om^2=&\,\der p-r\der x-rf'\der y,\\
\om^3=&\,\der r+r^2\frac{f^{(3)}}{f''}\der x-r^2\frac{f''{}^2-f'f^{(3)}}{f''}\der y,\\
\om^4=&\,\der x,\quad \om^5=\,\der y.
\end{aligned}
\]
Now, it is convenient to introduce \emph{new coordinates} $(X,Y,P,Q,R)$ on $M^5$ related to the coordinates $(x,y,z,p,r)$ via:
\[
\begin{aligned}
x=&\,-P+\frac{q f'}{f''},\\
y=&\,-\frac{q}{f''},\\
z=&\,Y-PX+q\frac{Xf'-f}{f''},\\
p=&\,X,\\
r=&\,\frac{1}{q-Q},
\end{aligned}
\]
where now, due to $p=X$, we have $f=f(X)$, $f'=f'(X)$ and
$f''=f''(X)$.  In these new coordinates the coframe
$(\om^1,\dots,\om^5)$ defining our para-CR structures reads:
\[
\begin{aligned}
\om^1=&\,\der Y-P\der X,\\
\om^2=&\,\frac{1}{q-Q}\,\Big(\,\der P-Q\der X\,\Big),\\
\om^3=&\,\frac{1}{(q-Q)^2}\,\Big(\,\der Q-\frac{f^{(3)}}{f''}\der P\,\Big),\\
\om^4=&\,-\der P+\der\, \big(q\frac{f'}{f''}\big),\\
\om^5=&\,-\frac{1}{f''}\,\Big(\der q-q\frac{f^{(3)}}{f''}\der X\,\Big).
\end{aligned}
\]
Now, a special para-CR transformation
\[
\begin{aligned}
&\om^1\to \om^1, \\& \om^2\to \big(q-Q\big)\,\om^2,\\
&\om^3\to \big(q-Q)^2\,\Big(\om^3+\frac{f^{(3)}}{(q-Q)f''}\om^2\Big),\\
&\om^4\to-\om^4-f'\om^5,\\
&\om^5\to -f''\,\om^5\end{aligned}
\]
as in~{\eqref{trump}}, brings this para-CR coframe to 
\[
\begin{aligned}
\om^1=&\,\der Y-P\der X,\\
\om^2=&\,\der P-Q\der X,\\
\om^3=&\,\der Q-Q\frac{f^{(3)}}{f''}\der X,\\
\om^4=&\,\der P-q\der X,\\
\om^5=&\,\der q-q\frac{f^{(3)}}{f''}\der X.
\end{aligned}
\]
Note the remarkable similarity of the 1-forms $(\om^2,\om^3)$ to the
1-forms $(\om^4,\om^5)$; they merely differ by the flip
$Q\leftrightarrow q$.

Now, let us consider two foliations of $M^5$ by two families of
hypersurfaces; $M^5$ is foliated by
\[
\mathrm{J^2}_{q_0}
=
\{M^5\in(X,Y,P,Q,q)\,:\, q=q_0=\mathrm{const}\}
\]
and by
\[
\mathrm{j^2}_{Q_0}=\{M^5\in(X,Y,P,Q,q)\,:\,
Q=Q_0=\mathrm{const}\}.\]
It follows from our calculations above that
every hypersurface $\mathrm{J^2}_{q_0}$ in the first family has a
structure of the space of second jets $\mathrm{J^2}$ coordinatized by
$(X,Y,P,Q)$ for the 3\textsuperscript{rd} order ODE
\begin{equation} 
Y'''=Y''\frac{f^{(3)}(X)}{f''(X)}\label{ee}
\end{equation}
for a function $Y=Y(X)$, with $Y'=P$, $Y''=Q$. Similarly, every
hypersurface $\mathrm{j^2}_{Q_0}$ in the second family has a structure
of the space of second jets $\mathrm{J^2}$ coordinatized by
$(X,Y,P,q)$ for the \emph{same} 3\textsuperscript{rd} order ODE $
Y'''=Y''\frac{f^{(3)}(X)}{f''(X)}$ for a function $Y=Y(X)$, with
$Y'=P$, $Y''=q$. Note that the passage from the first family of the
second jet spaces to the second family of the second jet spaces
corresponds to the flip $(\om^2,\om^3)\leftrightarrow (\om^4,\om^5)$
between the original coframe forms of the considered para-CR
structure~{\eqref{cre}}.

Since in our notation from Theorem~{\ref{th.c.1}} and
Proposition~{\ref{c.cor.har}} the ODE~{\eqref{ee}} has
$H=r\frac{f^{(3)}(x)}{f''(x)}$, then its Chern invariant
${\color{blue}\inc{Z}{}=H_{rrrr}\equiv 0}$. Thus, according to
Proposition~{\eqref{33}} each of the corresponding first jet spaces
$\mathrm{J^1}$ and $\mathrm{j^1}$, which curiously are both
parametrized by $(X,Y,P)$, has a natural projective contact structure.

To see this structure, we restrict to the case of $\mathrm{J^2}$; the case of $\mathrm{j^2}$ is the same, modulo the replacement $q\to Q$. Fortunately the ODE is easy to solve; its general solution is $$Y=c_1 f(X)+c_2 X+c3.$$ This general solution defines a 3-parameter family of curves
\[
\gamma(t;c_1,c_2,c_3)
=
(X(t),Y(t),P(t))=(t,c_1 f(t)+c_2 t+c_3,c_1f'(t)+c_2)
\]
in $\mathrm{J^1}$. Now we fix a frame $(e_1,e_2,e_3)=(\partial_Y,\partial_P,\partial_X+P\partial_Y)$ in $\mathrm{J^1}$, and consider tangent vectors $\dot{\gamma}(t)$ to each of these curves. Straightforward differentiation gives: 
$$\dot{\gamma}(t)=\sum_{i=1}^3\dot{\gamma}{}^ie_i=c_1f''(t)e_2+e_3.$$
Since the contact distribution $\mathcal{C}$ in $\mathrm{J^1}$ is given by
\[
\mathcal{C}
=
(\omega^1)^\perp=\Span(e_2,e_3)
\]
we see that our 3-parameter family of curves $\gamma(t;c_1,c_2,c_3)$ is always \emph{tangent} to $\mathcal{C}$. And now, writing the geodesic equations for the curves $\gamma(t)$ in the coframe $(e_1,e_2,e_3)$
\[
\frac{\der\dot{\gamma}{}^i}{\der t}
+
\sum_{j,k=1}^3\Gamma^i{}_{jk}\dot{\gamma}{}^j\dot{\gamma}{}^k
=
0,
\]
one can easilly see that our 3-parameter family of curves
$\gamma(t;c_1,c_2,c_3)$ satisfies these equations with a torsionless
connection $\nabla$, such that
$\nabla_{e_i}e_j=\sum_{k=1}^3\Gamma^k{}_{ji}e_k$, in which all the
coeffcients $\Gamma^k{}_{ij}=0$, except
$\Gamma^2{}_{23} = \Gamma^2_{32}=-\frac{f^{(3)}(X)}{2f''(X)}$.

Thus we have a 3-parameter family of curves $\gamma(t;c_1,c_2,c_3)$ in
$\mathrm{J^1}$, which are (a) tangent to the contact distribution
$\mathcal C$ and (b) are geodesics with respect to the torsionless
connection $\nabla$. This shows that $\mathrm{J^1}$ is equipped with a
\emph{contact projective structure}.

We thus have shown on an example, how a PDE
system~{\eqref{ss}}--{\eqref{ic}}, {\eqref{1deg}}--{\eqref{2ndeg}},
with ${\color{blue}2G_{ppp}+G_{pp}H_{rr}\equiv 0}$ defines two contact
equivalence classes of 3\textsuperscript{rd} odrer ODEs and a contact
projective structure on their space of first jets.

Finally, note that the quotient 3-manifolds on which the contact
projective structures associated with our para-CR structure resides
are just the quotients of the $M^5$ by the respective \emph{integrable
para-CR distributions} $D_1$ and $D_2$ in $M^5$.
\end{example}

%%%%%%%%%%%%%%%%%%%%%%%%%%%%%
\section{Appendix}\label{api}
%%%%%%%%%%%%%%%%%%%%%%%%%%%%%

It is instructive to show the result of Cartan's equivalence procedure
applied to the 1-forms~{\eqref{godwun}} or~{\eqref{godwn}} when we
have ${\color{blue} I^3 \equiv 0}$. We do it here for the 
1-forms~{\eqref{godwun}}.

For this, we need the system~{\eqref{sysend}} and its integrability
conditions, as in Theorem~{\ref{susend}}, adapted to
${\color{blue} I^3 \equiv 0}$. This restricted to
${\color{blue} I^3 \equiv 0}$ system reads:
%.../ala_godlinski_wunschman_I3eq0.nb
\begin{equation}
\begin{aligned}
\der \bar{\om}^1=&-\bar{\om}^1\dz\varpi_1+\bar{\om}^2\dz\bar{\om}^4,\\
\der \bar{\om}^2=&-\bar{\om}^1\dz\varpi_3+\bar{\om}^2\dz(\varpi_2-\tfrac12\varpi_1)+\bar{\om}^3\dz\bar{\om}^4,\\
\der \bar{\om}^3=&-\bar{\om}^2\dz\varpi_3+2\bar{\om}^3\dz\varpi_2+{\color{red}I^1}\,\bar{\om}^1\dz\bar{\om}^4,\\
\der \bar{\om}^4=&-\bar{\om}^1\dz\varpi_4-\bar{\om}^4\dz(\varpi_2+\tfrac12\varpi_1)-\bar{\om}^2\dz\bar{\om}^5,\\
\der \bar{\om}^5=&\,\bar{\om}^4\dz\varpi_4-2\bar{\om}^5\dz\varpi_2+{\color{green}I^2}\,\bar{\om}^1\dz\bar{\om}^2,
\end{aligned}\label{sysendu}
\end{equation}
\[
\begin{aligned}
\der{\color{red}I^1}=&\,I^1{}_{|1}\bar{\om}^1+I^1{}_{|2}\bar{\om}^2+I^1{}_{|3}\bar{\om}^3+I^1{}_{|4}\bar{\om}^4-\tfrac32 {\color{red}I^1}\varpi_1-3{\color{red}I^1}\varpi_2,\\
\der {\color{green}I^2}=&\,I^2{}_{|1}\bar{\om}^1+I^2{}_{|2}\bar{\om}^2+I^2{}_{|4}\bar{\om}^4+I^2{}_{|5}\bar{\om}^5-\tfrac32 {\color{green}I^2}\varpi_1+3{\color{green}I^2}\varpi_2.
\end{aligned}
\]

Integrability conditions of these equations imply an existence of a
1-form $\varpi_5$ such that:
\begin{equation}
\begin{aligned}
\der \varpi_1
=&\,
\bar{\om}^1\dz\varpi_5+\bar{\om}^2\dz\varpi_4-\bar{\om}^4\dz\varpi_3,\\
\der \varpi_2
=&\,-
\tfrac12\bar{\om}^2\dz\varpi_4-\tfrac12\bar{\om}^4\dz\varpi_3-\tfrac{1}{2}I^2{}_{|5}\bar{\om}^1\dz\bar{\om}^2+\\
&\ \ \ 
\tfrac{1}{2}I^1{}_{|3}
\bar{\om}^1\dz\bar{\om}^4-\bar{\om}^3\dz\bar{\om}^5,\\
\der \varpi_3
=&\,
\varpi_3\dz(\tfrac12\varpi_1+\varpi_2)+\tfrac12\bar{\om}^2\dz\varpi_5+\bar{\om}^3\dz\varpi_4+(I^1{}_{|23}+I^2{}_{|45})\bar{\om}^1\dz\bar{\om}^2+\\
&\ \ \
I^2{}_{|5}\bar{\om}^1\dz\bar{\om}^3-I^1{}_{|2}\bar{\om}^1\dz\bar{\om}^4+{\color{red}I^1}\bar{\om}^1\dz\bar{\om}^5-\tfrac{1}{2}I^1{}_{|3}\bar{\om}^2\dz\bar{\om}^4,\\
\der \varpi_4
=&\,
\varpi_4\dz(\tfrac12\varpi_1-\varpi_2)+\tfrac12\bar{\om}^4\dz\varpi_5+\bar{\om}^5\dz\varpi_3+\\
&\ \ \
I^2{}_{|4}\bar{\om}^1\dz\bar{\om}^2-{\color{green}I^2}\bar{\om}^1\dz\bar{\om}^3+\tfrac{1}{2}I^2{}_{|5}\bar{\om}^2\dz\bar{\om}^4-I^1{}_{|3}\bar{\om}^1\dz\bar{\om}^5,\label{sysendu1}\\
\der \varpi_5
=&\,
\varpi_5\dz\varpi_1+2\varpi_4\dz\varpi_3-I^2{}_{|5}\bar{\om}^1\dz\varpi_3-3I^1{}_{|3}\bar{\om}^1\dz\varpi_4+(I^2{}_{|15}+2I^2{}_{|44})\bar{\om}^1\dz\bar{\om}^2-\\
&\ \ \ 
4I^2{}_{|4}\bar{\om}^1\dz\bar{\om}^3+(I^1{}_{|31}-2I^1{}_{|22}-2I^1{}_{|234}-2I^2{}_{|445}) \bar{\om}^1\dz\bar{\om}^4-2(I^1{}_{|2}+I^1{}_{|34})\bar{\om}^1\dz\bar{\om}^5-\\
&\ \ \
{\color{green}I^2}\bar{\om}^2\dz\bar{\om}^3+(I^1{}_{|23}+I^2{}_{|45})\bar{\om}^2\dz\bar{\om}^4-I^1{}_{|3}\bar{\om}^2\dz\bar{\om}^5+I^2{}_{|5}\bar{\om}^3\dz\bar{\om}^4-{\color{red}I^1}\bar{\om}^4\dz\bar{\om}^5, 
\end{aligned}
\end{equation}
and as before, the coefficients ${\color{red}I^1}$ and
${\color{green}I^2}$ are, modulo a scale, the respective basic para-CR
relative invariants ${\color{red}\inc{A}{}}$ and
${\color{green}\inc{B}{}}$ from
Theorem~{\ref{the1}}:
\[
\begin{aligned}
{\color{red}I^1}\sim&~{\color{red}9D^2H_r-27DH_p-18H_rDH_r+18H_pH_r+4H_r^3+54H_z},\\
{\color{green}I^2}\sim&~{\color{green}40G_{ppp}^3-45G_{pp}G_{ppp}G_{pppp}+9G_{pp}^2G_{ppppp}}.
\end{aligned}
\]  

There is only one way of forcing the forms~{\eqref{godwun}}, with
1-forms $(\bar{\om}{}^1, \bar{\om}{}^2, \bar{\om}{}^3, \bar{\om}{}^4,
\bar{\om}{}^5)$ described by the
EDS~{\eqref{sysendu}}--{\eqref{sysendu1}}, to satisfy the
system~{\eqref{e.c.dtheta_10d}}. Such a requirement determines all
$\theta^i$s and $\Omega_\mu$s uniquely. Explicitely
\[
\theta^i=g^i{}_j\bar{\om}{}^j,
\quad\quad\mathrm{for\,\,all}
\quad\quad i=1,2,\dots,5,
\]
with the reduced matrix $g=(g^i{}_j)$ equal to
\[
g
=
\bma-\rho^2&0&0&0&0\\
f_2&\rho\mathrm{e}^\phi&0&0&0\\
-\frac{f_2^2}{2\rho^2}&-\frac{f_2\mathrm{e}^\phi}{\rho}&-\mathrm{e}^{2\phi}&0&0\\
\bar{f}{}_2&0&0&\rho\mathrm{e}^{-\phi}&0\\
-\frac{\bar{f}{}_2^2}{2\rho^2}&0&0&-\frac{\bar{f}{}_2\mathrm{e}^{-\phi}}{\rho}&\mathrm{e}^{-2\phi}
\ema,
\]
and the remaining forms $\Omega_1,\dots,\Omega_5$ are as follows:
\[
\begin{aligned}
\Omega_1=&\,-u_1\theta^1-\frac{\bar{f}_2}{\rho^2}\theta^2+\frac{f_2}{\rho^2}\theta^4+\varpi_1+\der\log(\rho^2), \\~&\\
\Omega_2=&\,\frac{2I^1{}_{|3}\rho^3\mathrm{e}^\phi-3f_2^2\bar{f}_2}{6\rho^6}\theta^1-\frac{f_2\bar{f}_2+\rho^4u_1}{2\rho^4}\theta^2-\frac{\bar{f}_2}{\rho^2}\theta^3-\frac{f_2^2}{2\rho^4}\theta^4-\frac{f_2}{\rho^2}(\tfrac12\varpi_1+\varpi_2)-\frac{\mathrm{e}^\phi}{\rho}\varpi_3+\frac{f_2}{\rho^2}\der\log(\frac{\rho\mathrm{e}^\phi}{f_2}),\\~&\\
\Omega_3=&\,-\frac{f_2\bar{f}_2+\rho^4u_1}{2\rho^4}\theta^1-\frac{\bar{f}_2}{\rho^2}\theta^2+\tfrac12\varpi_1-\varpi_2+\der\log(\rho\mathrm{e}^\phi), \\~\\
\Omega_4=&\,\frac{3f_2\bar{f}{}_2^2-2I^2{}_{|5}\rho^3\mathrm{e}^{-\phi}}{6\rho^6}\theta^1+\frac{\bar{f}{}_2^2}{2\rho^4}\theta^2+\frac{f_2\bar{f}{}_2-\rho^4u_1}{2\rho^4}\theta^4+\frac{f_2}{\rho^2}\theta^5+\frac{\bar{f}_2}{\rho^2}(\varpi_2-\tfrac12\varpi_1)-\frac{\mathrm{e}^{-\phi}}{\rho}\varpi_4+\frac{\bar{f}_2}{\rho^2}\der\log(\frac{\rho\mathrm{e}^{-\phi}}{\bar{f}_2}),
\end{aligned}
\]
\[
\begin{aligned}
\Omega_5=&\Big(\tfrac12 u_1^2+\frac{2I^1{}_{|23}+4I^2{}_{|45}}{3\rho^4}-\frac{I^2{}_{|5}f_2\mathrm{e}^{-\phi}+I^1{}_{|3}\bar{f}{}_2\mathrm{e}^\phi}{\rho^5}+\frac{f_2^2\bar{f}{}_2^2}{\rho^8}\Big)\theta^1+\Big(\frac{\bar{f}{}_2u_1}{\rho^2}-\frac{\mathrm{e}^{-\phi}I^2{}_{|5}}{3\rho^3}+\frac{f_2\bar{f}{}_2^2}{\rho^6}\Big)\theta^2+\frac{\bar{f}_2^2}{\rho^4}\theta^3+\\&\Big(\frac{f_2^2\bar{f}{}_2}{\rho^6}-\frac{f_2u_1}{\rho^2}-\frac{\mathrm{e}^{\phi}I^1{}_{|3}}{3\rho^3}\Big) \theta^4+\frac{f_2^2}{\rho^4}\theta^5-u_1\varpi_1+\frac{2f_2\bar{f}_2}{\rho^4}\varpi_2+\frac{2\bar{f}_2\mathrm{e}^\phi}{\rho^3}\varpi_3-\frac{2f_2\mathrm{e}^{-\phi}}{\rho^3}\varpi_4+\frac{1}{\rho^2}\varpi_5-\der u_1-\\&\frac{2u_1}{\rho}\der\rho+\frac{\bar{f}_2}{\rho^4}\der f_2-\frac{f_2}{\rho^4}\der\bar{f}_2-\frac{2f_2\bar{f}_2}{\rho^4}\der\phi.
\end{aligned}
\]
The resulting EDS \eqref{e.c.dtheta_10d} for these forms reads:
\[
\begin{aligned}
\der\hc{1} 
=&\,
\vc{1}\w\hc{1}+\hc{4}\w\hc{2},\nonumber \\
\der\hc{2}
=&\,
\vc{2}\w\hc{1}+\vc{3}\w\hc{2}+\hc{4}\w\hc{3},\nonumber \\
\der\hc{3}
=&\,
\vc{2}\w\hc{2}+(2\vc{3}-\vc{1})\w\hc{3}+\frac{\mathrm{e}^{\phi}}{3\rho^3}\,I^1{}_{|3}\,\hc{2}\w\hc{1}
-\big(\frac{\mathrm{e}^\phi}{\rho}\big)^3\,{\color{red}I^1}\,\hc{4}\w\hc{1}, \nonumber \\
\der\hc{4}
=&\,
\vc{4}\w\hc{1}+(\vc{1}-\vc{3})\w\hc{4}+\hc{5}\w\hc{2},\nonumber \\
\der\hc{5}
=&\,
\vc{4}\w\hc{4}
+(\vc{1}-2\vc{3})\w\hc{5}-\big(\frac{\mathrm{e}^{-\phi}}{\rho}\big)^3\,{\color{green}I^2}\,\hc{1}\w\hc{2}+\frac{\mathrm{e}^{-\phi}}{3\rho^3}\,I^2{}_{|5}\,\hc{1}\w\hc{4}, \nonumber \\
\der\vc{1}
=&\,
\vc{5}\w\hc{1}+\vc{4}\w\hc{2}-\vc{2}\w\hc{4},\nonumber \\
\der\vc{2}
=&\,
(\vc{3}-\vc{1})\w\vc{2}+\tfrac{1}{2}\vc{5}\w\hc{2}+\vc{4}\w\hc{3}
+\inc{A}{3}\,\hc{1}\w\hc{2}+\inc{A}{4}\hc{1}\w\hc{4}, \nonumber \\
\der\vc{3}
=&\,
\tfrac{1}{2}\vc{5}\w\hc{1}+\vc{4}\w\hc{2}+\hc{5}\w\hc{3}
-\frac{\mathrm{e}^{-\phi}}{3\rho^3}\,I^2{}_{|5}\,\hc{1}\w\hc{2}+\frac{\mathrm{e}^{\phi}}{3\rho^3}\,I^1{}_{|3}\,\hc{1}\w\hc{4},\\
\der\vc{4}
=&\,
\hc{5}\w\vc{2}+\vc{4}\w\vc{3}+\tfrac{1}{2}\vc{5}\w\hc{4}
+\inc{A}{6}\hc{1}\w\hc{2}-\inc{A}{3}\hc{1}\w\hc{4}, \nonumber \\
\der\vc{5}
=&\,
\vc{5}\w\vc{1}+2\vc{4}\w\vc{2}+\inc{C}{1}\hc{1}\w\hc{2}
+\inc{A}{8}\hc{1}\w\hc{4}. \nonumber
\end{aligned}
\]
Here
\[
\begin{aligned}
\inc{A}{3}=&\,\frac{\mathrm{e}^\phi\bar{f}{}_2I^1{}_{|3}}{3\rho^5}-\frac{\mathrm{e}^{-\phi}f_2I^2{}_{|5}}{3\rho^5}+\frac{I^1{}_{|23}}{3\rho^4}+\frac{I^2{}_{|45}}{3\rho^4},\\
\inc{A}{4}=&\,\frac{\mathrm{e}^\phi}{\rho^4}\Big(\frac{f_2I^1{}_{|3}}{3\rho}-\frac{\mathrm{e}^{2\phi}\bar{f}{}_{2}{\color{red}I^1}}{\rho}-\tfrac13\mathrm{e}^\phi\big(I^1{}_{|34}+3I^1{}_{|2}\big)\Big),\\
\inc{A}{6}=&\,\frac{\mathrm{e}^{-\phi}}{\rho^4}\Big(\frac{\bar{f}{}_2I^2{}_{|5}}{3\rho}-\frac{\mathrm{e}^{-2\phi}f_{2}{\color{green}I^2}}{\rho}+\tfrac13\mathrm{e}^{-\phi}\big(I^2{}_{|25}+4I^2{}_{|4}\big)\Big),\end{aligned}
\]
and we will not display $\inc{C}{1}$ and $\inc{A}{8}$ as not
important.

Since already here the symmetry ${\color{red}I^1} \leftrightarrow
{\color{green}I^2}$, corresponding to the change
$(\theta^2, \theta^3) \leftrightarrow (\theta^4,\theta^5)$, is clearly
visible, we skip writing down the analogous formulas for the 
1-forms~{\eqref{godwn}}.

\end{document}